\documentclass[11pt]{amsart}
\usepackage[pdftex,lmargin=1.4in, rmargin=1.4in,tmargin=1.25in,bmargin=1.25in, marginpar=2.9cm]{geometry}

\usepackage{amsmath}
\usepackage{amsfonts}
\usepackage{amssymb}
\usepackage{amsthm}

\usepackage[colorlinks=true,urlcolor=blue,linkcolor=blue,citecolor=blue,pagebackref,breaklinks]{hyperref}
\usepackage{cite}

\usepackage{mathrsfs}
\usepackage{bm}
\usepackage{dsfont}
\usepackage{upgreek}
\usepackage{xcolor}

\usepackage{graphicx}
\usepackage[all,cmtip]{xy}
\usepackage{enumerate}
\usepackage{blkarray}
\usepackage{array}
\usepackage{todonotes}

\renewcommand{\MR}[1]{ }

\usepackage{leftidx}
\usepackage{mathdots}
\usepackage{arydshln}
\usepackage{stmaryrd}

\newenvironment{psm}
{\left(\begin{smallmatrix}}
{\end{smallmatrix}\right)}

\newcommand{\IC}{\mathbb{C}}

\newcommand{\cmfield}{\cK}
\newcommand{\realfield}{\cmfield^+}
\newcommand{\CO}{\mathcal{O}}

\newcommand{\Heckealg}{\mathbb{T}}

\numberwithin{equation}{subsection}

\newtheorem{thm}{Theorem}[subsection]
\newtheorem{prop}[thm]{Proposition}

\newtheorem{coro}[thm]{Corollary}

\theoremstyle{remark}
\newtheorem{rmk}[thm]{Remark}

\newcommand{\mr}[1]{\mathrm{#1}}

\newcommand{\bA}{\mathbb{A}}

\newcommand{\bC}{\mathbb{C}}

\newcommand{\bQ}{\mathbb{Q}}
\newcommand{\bR}{\mathbb{R}}

\newcommand{\bV}{\mathbb{V}}

\newcommand{\bZ}{\mathbb{Z}}

\newcommand{\cD}{\mathcal{D}}

\newcommand{\cI}{\mathcal{I}}

\newcommand{\cK}{\mathcal{K}}

\newcommand{\cW}{\mathcal{W}}

\newcommand{\fs}{\mathfrak{s}}

\newcommand{\fB}{\mathfrak{B}}

\newcommand{\fH}{\mathfrak{H}}

\newcommand{\fQ}{\mathfrak{Q}}

\newcommand{\fS}{\mathfrak{S}}

\newcommand{\fU}{\mathfrak{U}}

\newcommand{\sB}{\mathscr{B}}

\newcommand{\sF}{\mathscr{F}}

\newcommand{\sS}{\mathscr{S}}

\DeclareMathAlphabet{\mathpzc}{OT1}{pzc}{m}{it}

\newcommand{\be}{\mathbf{e}}
\newcommand{\bfi}{\mathbf{i}}

\newcommand{\ue}{\underline{e}}
\newcommand{\uf}{\underline{f}}
\newcommand{\um}{\underline{m}}
\newcommand{\ul}{\underline{l}}

\newcommand{\utau}{\protect\underline{\tau}}
\newcommand{\unu}{\protect\underline{\nu}}

\DeclareMathOperator{\GL}{GL}

\DeclareMathOperator{\U}{U}

\DeclareMathOperator{\Sp}{Sp}

\DeclareMathOperator{\Mp}{Mp}

\newcommand{\hol}{\mathrm{hol}}

\newcommand{\sgn}{\mathrm{sgn}}

\newcommand{\Tr}{\mathrm{Tr}}

\newcommand{\vol}{\mathrm{vol}}
\newcommand{\Eis}{\mathrm{Eis}}

\DeclareMathOperator{\Ind}{Ind}

\DeclareMathOperator{\Lie}{Lie}

\DeclareMathOperator{\Sym}{Sym}

\newcommand{\lhra}{\ensuremath{\lhook\joinrel\longrightarrow}}
\newcommand{\lra}{\longrightarrow}
\newcommand{\ra}{\rightarrow}
\newcommand{\hra}{\hookrightarrow}

\newcommand{\ol}{\overline}

\newcommand{\wt}{\widetilde}

\newcommand{\llb}{\llbracket}
\newcommand{\rrb}{\rrbracket}

\newcommand{\bid}{\mathbf{1}}

\newdir{d}{{\subset}}   
\newcommand{\upl}{\mathrm{up}\text{-}\mathrm{left}}
\newcommand{\lowr}{\mathrm{low}\text{-}\mathrm{right}}

\makeatletter
\def\l@section{\@tocline{1}{0pt}{1pc}{}{}}
\def\l@subsection{\@tocline{2}{0pt}{1pc}{4.6em}{}}
\def\l@subsubsection{\@tocline{3}{0pt}{1pc}{7.6em}{}}
\renewcommand{\tocsection}[3]{%
  \indentlabel{\@ifnotempty{#2}{\makebox[2.3em][l]{%
    \ignorespaces#1 #2.\hfill}}}#3}
\renewcommand{\tocsubsection}[3]{%
  \indentlabel{\@ifnotempty{#2}{\hspace*{2.3em}\makebox[2.3em][l]{%
    \ignorespaces#1 #2.\hfill}}}#3}
\renewcommand{\tocsubsubsection}[3]{%
  \indentlabel{\@ifnotempty{#2}{\hspace*{4.6em}\makebox[3em][l]{%
    \ignorespaces#1 #2.\hfill}}}#3}
\makeatother

\title[Archimedean Zeta Integrals]{Archimedean Zeta Integrals for Unitary Groups}

\date{\today}

\author[E.\  Eischen]{E.\ Eischen$^\dag$} 
\thanks{$^\dag$Partially supported by NSF Grants DMS-1751281 and DMS-2302011.}
\author[Z.\ Liu]{Z.\ Liu$^\ddag$}
\thanks{$^\ddag$Partially supported by NSF Grant DMS-2001527}

\address{E. E. Eischen\\
Department of Mathematics\\
University of Oregon\\
Fenton Hall\\
Eugene, OR 97403-1222\\
USA}
\email{eeischen@uoregon.edu}

\address{Z. Liu\\
Department of Mathematics\\
South Hall, Room 6607\\
University of California\\
Santa Barbara, CA 93106-3080
}
\email{liuzheng072@gmail.com}

\begin{document}
\bibliographystyle{amsalpha} 

\maketitle
\vspace{-0.25in}

\begin{abstract}
We derive precise formulas for the archimedean Euler factors occurring in certain standard Langlands $L$-functions for unitary groups.  In the 1980s, Paul Garrett, as well as Ilya Piatetski-Shapiro and Stephen Rallis (independently of Garrett), discovered integral representations of automorphic $L$-functions that are Eulerian but, in contrast to the Rankin--Selberg and Langlands--Shahidi methods, do not require that the automorphic representations to which the $L$-functions are associated are globally generic.  Their approach, the {\it doubling method}, opened the door to a variety of applications that could not be handled by prior methods.  

For over three decades, though, the integrals occurring in the Euler factors at archimedean places for unitary groups eluded precise computation, except under particular simplifications (such as requiring certain representations to be one-dimensional, as Garrett did in the first major progress on this computation and only prior progress for general signatures).  We compute these integrals for holomorphic discrete series of general vector weights for unitary groups of any signature.  This has consequences not only for special values of $L$-functions in the archimedean setting, but also for $p$-adic $L$-functions, where the corresponding term had remained open.

\end{abstract}

\tableofcontents 
\numberwithin{equation}{subsection}

\section{Introduction}
In the 1980s, Paul Garrett, as well as Ilya Piatetski-Shapiro and Stephen Rallis (independently of Garrett), discovered integral representations of automorphic $L$-functions that are Eulerian but, in contrast to the Rankin--Selberg and Langlands--Shahidi methods, do not require that the automorphic representations to which the $L$-functions are associated are globally generic \cite{PSR, ga, cogdell, langlands, shahidi}.  Their approach, the {\it doubling method}, opened the door to a variety of applications that could not be handled by the Rankin--Selberg method, since not all automorphic representations of classical groups are globally generic.  Their work, together with later work of \cite{JSLi}, computed the  integrals occurring in each of the Euler factors produced by the doubling method at finite unramified primes.

For over three decades, though, the integrals occurring in the Euler factors at archimedean places eluded precise computation in the case of unitary groups, except in special cases that admitted particular simplifications (such as requiring certain representations to be one-dimensional, as Garrett did in the first major progress on this computation and only prior progress for general signatures).  The present paper computes them for holomorphic discrete series of general vector weights.  This has consequences not only for special values of $L$-functions in the archimedean setting, but also for $p$-adic $L$-functions, where the corresponding term had remained open except in special cases (as discussed in \cite[Section 4.5.4]{HELS}).

Applications of the doubling method include analytic, algebraic, and $p$-adic aspects of cuspidal automorphic $L$-functions.  For example, the doubling method was employed in the study of possible poles of automorphic $L$-functions in \cite{KR}.  Also, adapting Shimura's use of Rankin--Selberg integrals to prove the algebraicity of special values of (normalized) Rankin--Selberg convolutions in \cite{shimura76, shimura-hilbert}, Shimura and Harris employed doubling integrals in their proofs of the algebraicity of (normalized) critical values of automorphic $L$-functions in \cite{hasv, harriscrelle, shar}.  These proofs of algebraicity eventually led to constructions of $p$-adic $L$-functions.  In particular, in analogue with how Hida built on Shimura's work in the Rankin--Selberg case to construct $p$-adic $L$-functions associated to modular forms \cite{hi85}, several mathematicians, including each of the authors of the present paper, have constructed $p$-adic $L$-functions via the doubling construction. In the symplectic case, this was carried out in \cite{BS, liuJussieu}, while it was recently completed in the unitary case in \cite{HELS, ew}.

When these last two constructions of $p$-adic $L$-functions were completed, the precise form of their archimedean Euler factors remained open, except under restrictive assumptions.  For the case of symplectic groups, this problem was solved in \cite{liu-archimedean}. (With small adaptations, the method {\it loc.cit} also handles the unitary groups with signature $(a,a)$.) The general case of unitary groups, however, presents additional challenges.  In particular, the {\em Schr\"odinger model}, the tool employed in \cite{liu-archimedean}, lacks convenient, explicit formulas for unitary groups of signature $(a,b)$ with $a\neq b$.  In the present paper, we explain how to overcome these challenges and determine the precise form of the archimedean Euler factors appearing in \cite{HELS, ew} for unitary groups of arbitrary signature.

\subsection{Doubling method}
The doubling method (an instance of the {\it pullback methods} first discovered by Garrett in the 1970s, though years would pass before they appeared in the literature) takes as input several pieces of data meeting certain conditions, and it outputs an integral that has properties associated with a zeta function, namely a functional equation, a meromorphic continuation, and an Euler product.  One of the key steps in applying the doubling method involves choosing this data so that it is amenable to computation of these integrals.  In addition, one must carefully choose data that result in $L$-functions suited to one's desired applications.

Before proceeding, we briefly recall the construction of the doubling integral.  (Later in the paper we give more details, but at the moment, we give only the information necessary for clarity in our introduction to the problem.  See also \cite{cogdell} for an introduction to the doubling method.)  Consider a Hecke character $\chi$ of a CM field $\cmfield$ with maximal totally real subfield $\realfield$, and consider a unitary group $G$ that preserves a Hermitian pairing on a vector space of defined over $\cmfield$ of dimension $n$.  Then $G\times G$ can be identified with a subgroup of a unitary group $H$ of signature $(n,n)$ (at each archimedean place). For the doubling method, we also need a particular Siegel parabolic subgroup $P$ of $H$.  (So the Levi subgroup $M_P$ is isomorphic to $\GL(n)_{\cmfield}$.)  
A Siegel Eisenstein series on $H$ associated to $s\in \IC$ and a section
$f(s, \chi)=\otimes' f_v(s,\chi)\in \Ind_{P\left(\bA_{\realfield}\right)}^{H\left(\bA_{\realfield}\right)}\left(\chi\cdot\left|\cdot\right|^{s}\right)\cong\otimes_v\Ind_{P\left(\realfield_v\right)}^{H(\realfield_v)}\left(\chi_v\cdot\left|\cdot\right|_v^{s}\right)$
(with the product over all places of $\realfield$ and $\chi\cdot \left|\cdot\right|^{s}$ viewed as a character on $P$ via composition with the projection from $P$ to $M_P\cong\GL(n)_{/\cmfield}$ and the determinant map on $\GL(n)_{/\cmfield}$) is the $\IC$-valued function
\begin{align}\label{siegelseries}
E\left(g, f(s, \chi)\right) = \sum_{\gamma\in P\left(\realfield\right)\backslash H\left(\realfield\right)}f(s, \chi)\left(\gamma g\right).
\end{align}
Let $\pi$ be a cuspidal automorphic representation of $G$ and $\tilde{\pi}$ be its contragredient.  Let $\varphi\in\pi$ and $\tilde{\varphi}\in\tilde{\pi}$.
The global doubling integral is then
\begin{align}\label{doublequ}
Z\left(f(s,\chi), \varphi, \tilde{\varphi}\right) = \int_{[G\times G]}\varphi(g)\tilde{\varphi}(h)E\left(\left(g, h\right), f(s, \chi)\right)\chi^{-1}\left(\det(h)\right)dg dh,
\end{align}
with $E(\cdot, f(s, \chi))$ an Eisenstein series as in Equation \eqref{siegelseries}.  As proved in \cite{PSR}, when these data are chosen so that we have factorizations into restricted tensor products $\varphi = \otimes'_v\varphi_v,$ $\tilde{\varphi} = \otimes'_v\tilde{\varphi}_v$, and $\otimes'_vf(s, \chi_v)$, this integral unfolds as a product of local integrals $Z\left(f(s, \chi), \varphi, \tilde{\varphi} \right)=\prod_v Z_v\left(f_v(s, \chi), \varphi_v, \tilde{\varphi}_v \right)$ with
\begin{equation}\label{Zv}
     Z_v\left(f_v(s, \chi), \varphi_v, \tilde{\varphi}_v \right)=\int_{G(\realfield_v)} f_v(s,\chi)\big((g,1)\big)\left<\pi_v(g)\varphi_v,\wt{\varphi}_v\right>\,dg, 
\end{equation}
where $\langle, \rangle$ denotes a suitably normalized $G$-invariant pairing.  For finite places $v$ where everything is unramified, the integral in Equation \eqref{Zv} was computed in the split case in \cite{PSR} and in the inert case in \cite{JSLi} and shown to equal $L_v\left(s+\frac{1}{2},\pi_v\times\chi_v\right)$ up to a product of explicit $L$-factors for characters. (More recently, Yamana proved that for all finite places, the greatest common denominator of the local zeta integrals is the correct definition of local $L$-factors \cite[Theorem 5.2]{Yamana}.) So the behavior of the global $L$-function is determined by the global integral \eqref{doublequ} and thus by the behavior of the cusp forms and the pullback of the Eisenstein series to $G\times G$ if one has a good understanding of the local integral \eqref{Zv}.

The constructions of $p$-adic $L$-functions via the doubling method, e.g. in \cite{liuJussieu, HELS}, depend on carefully choosing the local data. For $p$-adic $L$-functions in the setting of unitary groups, the local zeta integrals at all the finite places are calculated in \cite{HELS} (for the local data chosen there).  On the other hand, precise values at the archimedean places were only known in certain cases.  In the absence of precise formulations of the archimedean Euler factors, one has only a poorly understood normalizing factor in the comparison of the $p$-adic $L$-function with the corresponding $\IC$-valued $L$-functions, whose critical values the $p$-adic $L$-function interpolates.

\subsection{Relationship with prior results}
Prior to the present paper, some special cases had been addressed.  In \cite{ga06}, Garrett showed that the archimedean zeta integrals are algebraic up to a predictable power of the transcendental number $\pi$.  At each archimedean place where the extreme $K$-type of the cuspidal automorphic representation $\pi$ is one-dimensional, a precise computation of the archimedean zeta integrals is given in \cite{sh, shar}.  More generally, at each archimedean place at which at least one of the two factors of the extreme $K$-type of $\pi$ is one-dimensional, Garrett computed the archimedean integrals precisely for a certain choice of sections.  The general case, where neither factor was required to be a scalar, though, remained open until the present paper.

Key inspiration for the approach of the present paper comes from the theta correspondence. The computation of archimedean zeta integrals from the point of view of the theta correspondence has been studied in \cite{liu21, liun1}, where B. Lin and D. Liu consider unitary groups of signature $(n,1)$ and the evaluation of the archimedean zeta integrals at the center. The main interest of the computation {\it loc. cit.} is about non-holomorphic discrete series. Like \cite{liu21, liun1}, we will employ the {\it Fock} models, perhaps more familiar in certain areas of harmonic analysis and mathematical physics than in algebraic number theory.  

In the symplectic case in \cite{liu-archimedean}, the Schr\"odinger model is employed.  As noted above, though, the Schr\"odinger model alone lacks convenient formulas for unitary groups of signature $(a,b)$ with $a\neq b$.  Consequently, in the present paper, we must employ both the Fock and Schr\"odinger models, which are related through the {\em Bargmann transform}.
 
\subsection{Main results and overview of paper}

In \S \ref{archimedean-section}, after introducing notation and conventions, we introduce the archimedean zeta integrals  that we will compute.  In particular, we introduce the doubling integral and our precise choices of archimedean data.  These choices are based on those from \cite{HELS}.  In particular, we work with  holomorphic discrete series.

Section \ref{theta-section} introduces the Schr\"odinger and Fock models, which provide convenient models for the Weil representation and are related through the Bargmann transform.  
By using the Schr\"odinger model (\S \ref{section:schrodinger}), we reduce the computation of the archimedean zeta integrals to that of certain matrix coefficients of the Weil representation. Unlike in the case of symplectic groups or unitary groups of signature $(a,a)$, the Schr\"odinger model (\S \ref{section:schrodinger}) lacks convenient, explicit formulas (and pluri-harmonic polynomials) to compute those matrix coefficients. We use the Bargmann transform to carry out the computation with the Fock model.

Section \ref{section:threespaces} further simplifies the integrals. The main trick is to consider three different one-dimensional spaces inside the subspace of the joint harmonic polynomials which is isomorphic to the tensor product of an irreducible representation of $\U(a)\times\U(b)\times\U(k)$ and its dual representation.

In \S \ref{computation-section}, we obtain our main result, Theorem \ref{thm:Zright}.  As a consequence, we see that the normalization factor we compute has the form predicted by Coates and Perrin-Riou in \cite{coates, CPR}.    In particular, concluding the calculation of the doubling integral at a critical point for a unitary group at an archimedean place, our main result is:

\begin{thm}[Theorem \ref{thm:Zright}]\label{thm:Z}
Let $\sigma$ be an archimedean place of $\cK^+$ where the unitary group $G(\cK^+_\sigma)$ has signature $(a,b)$ (with $n=a+b$). Denote by $\cD_{(\utau;\unu)}$ (resp. $\cD^*_{(\utau;\unu)}$) the holomorphic discrete series of $G(\cK^+_\sigma)$ of weight $(\utau;\unu)=(\tau_1,\dots,\tau_a;\nu_1,\dots\nu_b)$ (resp. the contragredient representation of $\cD_{(\utau;\unu)}$), and denote by $\chi_{\mr{ac}}$ the anticyclotomic character of $\bC^\times$ sending $x$ to $\frac{x}{|x\ol{x}|^{1/2}}$. Suppose that the integer $k\geq n$ satisfies the Condition \eqref{kcondition}. 

Then for the archimedean section $f_{k,(\utau;\unu)}(s,\chi^r_{\mr{ac}})$ defined in Equation \eqref{eq:fktv}, and $v_{(\utau;\unu)}\in \cD_{(\utau;\unu)}$ (resp. $v^*_{(\utau;\unu)}\in \cD^*_{(\utau;\unu)}$), the highest weight vector inside the lowest $K$-type of $\cD_{(\utau;\unu)}$ (resp. the dual vector of $v_{(\utau;\unu)}$), 
{\small
\begin{align*}
   &\left.Z_\sigma\left(f_{k,(\utau;\unu)}\left(s,\chi^r_{\mr{ac}}\right),v^*_{(\utau;\unu)},v_{(\utau;\unu)}\right)\right|_{s=\frac{k-n}{2}}\\
   =&\,\frac{2^{ab-
   \frac{n}{2}}\pi^{ab}\,(2\pi i)^{-\sum\tau_j-\sum\nu^*_j+\frac{a(k+r)}{2}+\frac{b(k-r)}{2}}}{\dim\left(\GL(a),\utau\right)\dim\left(\GL(b),\unu\right)}\frac{\prod\limits_{j=1}^a\Gamma\left(\tau_j-j+\frac{k-r}{2}-b+1\right)\prod\limits_{j=1}^b\Gamma\left(\nu^*_j-j+\frac{k+r}{2}-a+1\right)}{\prod\limits_{j=1}^n\Gamma(k-j+1)}\\
   =&\,\frac{ 2^{\frac{n^2}{2}-n}i^{-\frac{n^2}{2}+\frac{n}{2}-ab}\,(-1)^{nr}}{{\dim\left(\GL(a),\utau\right)\dim\left(\GL(b),\unu\right)}}\frac{E_\sigma\left(\frac{k-n+1}{2},\cD^*_{(\utau;\unu)}\times\chi^r_{\mr{ac}}\right)}{2^{n(n-1)}\pi^{\frac{n(n-1)}{2}}(-2\pi i)^{-nk}\prod_{j=1}^n\Gamma(k+1-j)}.
\end{align*}
}

Here, the term $E_\sigma\left(s,\cD^*_{(\utau;\unu)}\times\chi^r_{\mr{ac}}\right)$ denotes the modified archimedean Euler factor defined in \S \ref{sec:modiE}.
\end{thm}

Condition \eqref{kcondition} is for the point $s=\pm\frac{k-n}{2}+\frac{1}{2}$ to be critical for the $L$-factor $L_\sigma\left(s,\cD^*_{(\utau;\unu)}\times\chi^r_{\mr{ac}}\right)$. The factor 
\begin{align*}
2^{n(n-1)}\pi^{\frac{n(n-1)}{2}}(-2\pi i)^{-nk}\prod_{j=1}^n\Gamma(k+1-j)
\end{align*}
also arises in  \cite[Equation (15)]{apptoSHL} in the process of constructing  a family of Eisenstein series that can be $p$-adically interpolated, where (like in the present paper) it is the archimedean normalization factor for Siegel Eisenstein series at $s=\frac{k-n}{2}$ on $\U(n,n)$.  
Combining Theorem~\ref{thm:Z} and the functional equation for doubling local zeta integrals, one deduces formulas for 
\[
   \left.Z_\sigma\left(f_{k,(\utau;\unu)}\left(s,\chi^r_{\mr{ac}}\right),v^*_{(\utau;\unu)},v_{(\utau;\unu)}\right)\right|_{s=\frac{n-k}{2}}
\] 
in terms of modified Euler factor $E_\sigma\left(\frac{n-k+1}{2},\cD^*_{(\utau;\unu)}\times\chi^r_{\mr{ac}}\right)$ for integers $k$ satisfying Condition \eqref{kcondition}. See Theorem~\ref{thm:left} for details.

Applying our main theorem to the $p$-adic measures from \cite[Main Theorem 9.2.2]{HELS}, we obtain Corollary \ref{coro:apptoHELS} below.  Before stating it, we make a few remarks about its contents.  Because it would take considerable space to fully explain each of the technical conditions of \cite[Main Theorem 9.2.2]{HELS} and they are unnecessary for developing the main results for the present paper, we just briefly highlight them here, to give the reader a sense of the setting in which Corollary \ref{coro:apptoHELS} holds.   The representation $\pi$ is assumed to be of type $(\kappa, K_{r_1})$ in the sense of \cite[Section 6.5.1]{HELS}, where $K_{r_1}$ is a certain open compact subgroup and $r_1$ is a parameter tied to the level.  (This condition in particular implies that all the archimedean components of $\pi$ are isomorphic to holomorphic discrete series.)  $\Heckealg=\Heckealg_\pi$ denotes the corresponding connected component of the ordinary Hecke algebra, and $\lambda_{\pi}$ is a character of $\Heckealg$ as in \cite[Section 6.6.8]{HELS}.  The main hypotheses on $\pi$ are that it satisfies a global multiplicity one hypothesis (\cite[Hypothesis 7.3.3]{HELS}), a Gorenstein condition (\cite[Hypothesis 7.3.2]{HELS}), and a minimality hypothesis (\cite[Proposition 7.3.5]{HELS}).  The cuspidal automorphic forms $\varphi$, $\varphi^\flat$ are chosen to lie in certain lattices.  

We define $X_p$ to be the maximal abelian extension of $\cmfield$ unramified away from $p$, and $\Lambda_{X_p}$ denotes an Iwasawa algebra $\CO\llb X_p\rrb$ for a sufficiently large $p$-adic ring $\CO$.  We also take $S$ to be a finite set of primes including all ramified places and all places dividing $p\infty$, and we denote by $L^S$ the (product of the Euler factors of the) standard Langlands $L$-function away from $S$.

For the automorphic representation $\pi$ as above and an algebraic unitary Hecke character $\chi_u$ of $\cK^\times\backslash\bA^\times_{\cK}$, we set
\begin{align*}
E_\infty&\left(s, \pi, \chi_u\right):=\prod_{\sigma:\realfield\hookrightarrow\bR}E_\sigma\left(s,\pi^*_\sigma\times\chi_{u, \sigma}\right).
\end{align*}
We also define the constant 
\[
   A\left(\pi_\infty,\chi_{u,\infty}\right):=\prod_{\sigma:\realfield\hookrightarrow\IC}\frac{ 2^{\frac{n^2}{2}-n}i^{-\frac{n^2}{2}+\frac{n}{2}-a_\sigma b_\sigma}\,\chi_{u,\sigma}(-1)^{n}}{{\dim\left(\GL\left(a_\sigma\right),\utau_\sigma\right)\dim\left(\GL\left(b_\sigma\right),\unu_\sigma\right)}},
\]
where $(a_\sigma,b_\sigma)$ is the signature of the unitary group at $\sigma$, and $\left(\utau_\sigma;\unu_\sigma\right)$ is the weight of the holomorphic discrete series $\pi_\sigma$.

In addition, to emphasize the connection with the Euler factors predicted by Coates in \cite{coates}, we follow his conventions, writing $E_p$ for the product of the factors at primes $v$ dividing $p$ in place of the notation employed in \cite[Equation (86)]{HELS}.

\begin{coro}\label{coro:apptoHELS}
Let $\varphi$ and $\varphi^\flat$ be cuspidal automorphic forms from $\pi$ and its contragredient, respectively, 
all meeting the conditions of \cite[Main Theorem 9.2.2]{HELS}. There is a unique element $L\left(\Eis, \varphi\otimes\varphi^\flat\right)\in \Lambda_{X_p}\hat{\otimes}\Heckealg$ such that for any algebraic Hecke character $\chi = ||\cdot||^{-\frac{k-n}{2}}\chi_u$ of $X_p$, with $\chi_u$ unitary and $k\geq n$ a positive integer, the image of $L\left(\Eis, \varphi\otimes\varphi^\flat\right)$ under the map on $\Lambda_{X_p}\hat{\otimes}\Heckealg$ induced by $\chi\otimes\lambda_{\pi}$ is
\begin{align*}
\Omega_{\varphi, \varphi^\flat}  I_SA(\pi_\infty,\chi_{u,\infty})
 L^S\left(\frac{k-n+1}{2}, \pi, \chi_u\right)E_p\left(\frac{k-n+1}{2}, \pi, \chi_u\right)E_\infty\left(\frac{k-n+1}{2}, \pi, \chi_u\right),
\end{align*}
with $\Omega_{\varphi, \varphi^\flat}$ a period associated to $\varphi$ and $\varphi^\flat$ (more precisely, the Petersson pairing of these two forms, normalized by a volume factor) and $I_S$ a finite product (of certain constant volume factors and Euler factors at finite primes in $S$ of Hecke $L$-functions associated to $\chi$). 
\end{coro}
In the formulation of \cite[Main Theorem 9.2.2]{HELS}, the period $\Omega_{\varphi, \varphi^\flat}$ is instead expressed as a product of three terms that are individually dependent on $\chi$ but whose product ($\Omega_{\varphi, \varphi^\flat}$) is independent of $\chi$.  (The notation $\Omega_{\varphi, \varphi^\flat}$ is not used in that paper, but we use it here to emphasize the dependence only on $\varphi$ and $\varphi^\flat$.  The term $\Omega_{\varphi, \varphi^\flat}$ is the product of the first two terms, a volume factor and a pairing, from the expression obtained in \cite[Corollary 9.2.1]{HELS}.  As explained in \cite[Remark 4.1.6]{HELS}, this pairing is the Petersson pairing of $\varphi$ and $\varphi^\flat$.)  For the reader expecting to see a Gauss sum, please note that they arise here not in $\Omega_{\varphi, \varphi^\flat}$ but instead in the $\varepsilon$-factors in $E_p\left(\frac{k-n+1}{2}, \pi, \chi_u\right)$.  (Similarly, Gauss sums arise in the Euler factors at $p$ in the analogous construction for symplectic groups in \cite{liuJussieu} whose archimedean factors were computed in \cite{liu-archimedean}.)
The notation $\Eis$ in $L\left(\Eis, \varphi\otimes\varphi^\flat\right)$ refers to the {\it Eisenstein measure} from \cite{apptoSHL, apptoSHLvv} used to construct the $p$-adic $L$-functions in \cite{HELS}.

\subsection{Acknowledgements}
We thank Paul Garrett, Michael Harris, and Dongwen Liu for helpful conversations related to this problem.  We are also grateful to the referee for carefully reading the paper and providing helpful feedback.

\section{The archimedean zeta integrals for studying critical $L$-values}\label{archimedean-section}
\subsection{The unitary group $\U(a,b)$}\label{notation-section}
We begin by introducing some notation and conventions for unitary groups.
Let $a,b$ be two non-negative integers with $a+b=n$, and
$$\U(a, b)=\left\{g\in\GL(n,\bC)\middle|\,\ltrans{\ol{g}}\begin{pmatrix}\bid_a&\\&-\bid_b\end{pmatrix}g=\begin{pmatrix}\bid_a&\\&-\bid_b\end{pmatrix}\right\},$$
the unitary group of signature $(a, b)$.
We fix the Haar measure on $\U(a,b)$ as the product measure of the maximal compact subgroup $\U(a)\times\U(b)$ and the symmetric domain $\fB_{a,b}=\left\{z\in M_{a,b}(\bC)\middle|\,\bid_a-z\ltrans{\ol{z}}>0\right\}$, where the measure on $\U(a)\times\U(b)$ is the Haar measure with total volume $1$, and the measure on $\fB_{a,b}$ is 
$$\det(1-z\ltrans{\ol{z}})^{-n}\prod\limits_{\substack{1\leq i\leq a\\1\leq j\leq b}}|dz_{ij}d\ol{z}_{ij}|.$$
This Haar measure agrees with the standard Haar measure (see, e.g., \cite[\S7]{HC-HI})  with respect to the symmetric form $(X,Y)=\Tr\ltrans{\ol{X}}Y$ on $\Lie\U(a,b)$.
\subsection{Archimedean $L$-factors for holomorphic discrete series}

We briefly recall the definition of the local $L$-factors for the holomorphic discrete series representation of $\U(a,b)$ with which we work. Let $\cD_{(\utau;\unu)}$ be the holomorphic discrete series of weight $(\utau;\unu)=(\tau_{1},\ldots,\tau_{a};\nu_{1},\ldots,\nu_{b})$. We denote by $\cD^*_{(\utau;\unu)}$ the contragredient representation of $\cD_{(\utau;\unu)}$ which is the anti-holomorphic discrete series of weight $(\utau^*;\unu^*)$. We write $\utau^*=(\tau^*_1,\dots,\tau^*_a)=(-\tau_a,\dots,-\tau_1)$ and $\unu^*=(\nu^*_1,\dots,\nu^*_b)=(-\nu_b,\dots,-\nu_1)$. 

The half sum of compacts roots (resp. non-compact roots) of $\U(a,b)$ is
\begin{align*}
   \rho_{\mr{c}}&=\bigg(\underbrace{\frac{a-1}{2},\dots,\frac{-a+1}{2}}_a;\underbrace{\frac{b-1}{2},\dots,\frac{-b+1}{2}}_b\bigg) &(\text{resp. }\rho_{\mr{nc}}&=\bigg(\underbrace{\frac{b}{2},\dots,\frac{b}{2}}_a;\underbrace{-\frac{a}{2},\dots,-\frac{a}{2}}_b\bigg)),
\end{align*} 
and the Harish-Chandra parameter of $\cD_{(\utau;\unu)}$ equals
\begin{equation}\label{eq:HCp}
\begin{aligned}
   &(\utau;\unu)+(\rho_{\mr{c}}-\rho_{\mr{nc}})\\
   =&\,\left(\tau_1+\frac{a-b-1}{2},\dots,\tau_a+\frac{-a-b+1}{2};\nu_1+\frac{a+b-1}{2},\dots,\nu_b+\frac{a-b+1}{2}\right).
\end{aligned}
\end{equation}

Let $\chi_{\mr{ac}}$ be the unitary character of $\bC^\times$ which takes the value $\frac{x}{(x\ol{x})^{1/2}}$ at $x\in\bC$, and $r$ be an integer. The local $L$-factor for $\cD^*_{(\utau;\unu)}\times\chi^r_{\mr{ac}}$ is
\begin{align*}
   &L\left(s,\cD^*_{(\utau;\unu)}\times\chi^r_{\mr{ac}}\right)\\
   =&\,\prod_{j=1}^a\Gamma_{\bC}\left(s+\left|\tau_j-\frac{r}{2}+\frac{a-b+1}{2}-j\right|\right)\prod_{j=1}^b\Gamma_{\bC}\left(s+\left|\nu_j-\frac{r}{2}+\frac{a+b+1}{2}-j\right|\right),
\end{align*}
where $\Gamma_{\bC}(s)=2(2\pi)^{-s}\Gamma(s)$. Since we will focus on one real place, we omit the subscript indicating the place from the $L$-factor $L$, the modified Euler factor $E$ and the local doubling zeta integral.

The condition for $s=s_0$ to be a critical point for $L\left(s,\cD^*_{(\utau;\unu)}\times\chi^r_{\mr{ac}}\right)$ is that $$s_0+\frac{r+n-1}{2}\in\bZ,$$ 
and
\begin{align*}
   -\left|\tau_j-\frac{r}{2}+\frac{a-b+1}{2}-j\right|+1&\leq s_0\leq \left|\tau_j-\frac{r}{2}+\frac{a-b+1}{2}-j\right|, &\text{for all }1\leq j\leq a,\\
   -\left|\nu_j-\frac{r}{2}+\frac{a+b+1}{2}-j\right|+1&\leq s_0\leq \left|\nu_j-\frac{r}{2}+\frac{a+b+1}{2}-j\right|, &\text{for all }1\leq j\leq b.
\end{align*}

\subsection{Modified archimedean Euler factor for $p$-adic interpolation}\label{sec:modiE} Following the conventions of \cite{coates}, the modified archimedean Euler factor for $\cD^*_{(\utau;\unu)}\times\chi^r_{\mr{ac}}$ is 
\begin{equation}\label{eq:Esigma}
\begin{aligned}
   E\left(s,\cD^*_{(\utau;\unu)}\times\chi^r_{\mr{ac}}\right)=&\prod^a_{j=1}e^{-\frac{\pi i}{2}\left(s+\left|\tau_j-\frac{r}{2}+\frac{a-b+1}{2}-j\right|\right)}\,\Gamma_{\bC}\left(s+\left|\tau_j-\frac{r}{2}+\frac{a-b+1}{2}-j\right|\right)\\
   &\times\prod_{j=1}^b e^{-\frac{\pi i}{2}\left(s+\left|\nu_j-\frac{r}{2}+\frac{a+b+1}{2}-j\right|\right)}\,\Gamma_{\bC}\left(s+\left|\nu_j-\frac{r}{2}+\frac{a+b+1}{2}-j\right|\right).
\end{aligned}
\end{equation}
(In {\it loc.cit}, the motives are assumed to be defined over $\bQ$, and a subscript $\infty$ is used to denote the archimedean place.)

For an integer $k$ satisfying the later defined condition \eqref{kcondition}, expanding the right hand side of \eqref{eq:Esigma} at  $s=\frac{k-n+1}{2}$, we have
\begin{align*}
    E\left(\frac{k-n+1}{2},\cD^*_{(\utau;\unu)}\times\chi^r_{\mr{ac}}\right)=&\,2^{-n}(2\pi i)^{-\sum\tau_j-\sum\nu^*_j-a\frac{k-r}{2}-b\frac{k+r}{2}+\frac{a(a-1)+b(b-1)}{2}+2ab}\\
   &\hspace{-4em}\times\prod_{j=1}^a\Gamma\left(\tau_j+1-j+\frac{k-r}{2}-b\right)\prod_{j=1}^b\Gamma\left(\nu^*_j+1-j+\frac{k+r}{2}-a\right).
\end{align*}
Let $\gamma\left(s,\cD^*_{(\utau;\unu)}\times\chi^r_{\mr{ac}}\right)$ be the gamma factor as defined in \cite{LapidRallis}. Then it follows from the properties of the gamma factor and the definition of the modified Euler factor that
\begin{equation}\label{eq:Eleft}
\begin{aligned}
   E\left(\frac{n-k+1}{2},\cD^*_{(\utau;\unu)}\times\chi^r_{\mr{ac}}\right)=&\,(-1)^{\sum\tau_j+\sum\nu_j+a\frac{k+r}{2}+b\frac{k-r}{2}}i^{-a^2-b^2}\\
   &\hspace{-6em}\times\gamma\left(\frac{k-n+1}{2},\cD^*_{(\utau;\unu)}\times\chi^r_{\mr{ac}}\right)\,E\left(\frac{k-n+1}{2},\cD^*_{(\utau;\unu)}\times\chi^r_{\mr{ac}}\right).
\end{aligned}
\end{equation}

\subsection{The archimedean sections used for $p$-adic $L$-functions}
We now introduce the choices of archimedean sections used in the doubling integrals for the construction of the $p$-adic $L$-functions in \cite{HELS}.
Let $J_{2n}=\begin{pmatrix}&i\cdot\bid_n\\ -i\cdot\bid_n&\end{pmatrix}$. Viewing $J_{2n}$ as an Hermitian form on $\bC^{2n}$, the associated unitary group 
$$\U(J_{2n})=\left\{h\in\GL(2n,\bC):\,\ltrans{\ol{h}}J_{2n}h=J_{2n}\right\}$$
is isomorphic to $\U(n,n)$ and is quasi-split. It contains a Siegel parabolic subgroup
$$
   Q_{\U(J_{2n})}=\left\{h=\begin{pmatrix}A&B\\&\ltrans{\ol{A}}^{-1}  \end{pmatrix}:\,A\in\GL(n,\bC),\, B\in M_{n,n}(\bC),\, \ltrans{(A^{-1}B)}=\ol{A^{-1}B} \right\}.
$$

Let $I\left(s,\chi^r_{\mr{ac}}\right)=\Ind_{Q_{\U(J_{2n})}}^{\U(J_{2n})}\chi^r_{\mr{ac}}|\cdot|^s_{\bC}$ be the (normalized) induction of the character $\begin{pmatrix}A&B\\&\ltrans{\ol{A}}^{-1}  \end{pmatrix}\mapsto \chi^r_{\mr{ac}}(\det A)|\det A|^{s}_{\bC}$ (where $|x|_{\bC}=x\ol{x}$ for $x\in\bC$). Inside the degenerate principal series $I\left(s,\chi^r_{\mr{ac}}\right)$, for an integer $k$ of the same the parity as $r$, we have the classical section $f_k(s,\chi^r_{\mr{ac}})$ (built from the canonical automorphy factor) , whose values on $h =\begin{pmatrix}A&B\\C&D\end{pmatrix}$, with $A, B, C, D$ $n\times n$ matrices, are given by
\begin{equation}\label{eq:fk}
   f_k(s,\chi^r_{\mr{ac}})\left(h\right)=(\det h)^{\frac{r+k}{2}}\det(Ci+D)^{-k}|\det(Ci+D)|^{-s+\frac{k}{2}-\frac{n}{2}}_{\bC}.
\end{equation}

In \cite{HELS}, the section $f_{k,(\utau;\unu)}(s,\chi^r_{\mr{ac}})\in I(s,\chi^r_{\mr{ac}})$ chosen for $s=\frac{k-n}{2}$ with $k\geq n$ and $\cD^*_{(\utau;\unu)}$ twisted by $\chi^r_{\mr{ac}}$ is obtained by applying the Lie algebra operators to the classical section $f_k(s,\chi^r_{\mr{ac}})$. 

Let
\begin{equation}\label{eq:s}
   \fs=\frac{1}{\sqrt{2}}\begin{pmatrix}\bid_n&-i\cdot\bid_n\\\bid_n&i\cdot\bid_n\end{pmatrix},
\end{equation}
the conjugation of which gives an isomorphism between $\U(n,n)$ and $\U(J_{2n})$. We consider the Lie algebras $\big(\Lie\U(J_{2n})\big)\otimes_{\bR}\bC$ and $\big(\Lie\U(n,n)\big)\otimes_{\bR}\bC$. In order to distinguish from the $i$ in the first factor of the tensor product, we write $\bfi$ for a fixed square root of $-1$ in the second factor. The weight raising  operators in $\left(\Lie \U(n,n)\right)\otimes_{\bR}\bC$ consist of
\begin{align*}
   \mu^+_{\U(n,n),X}&=\frac{1}{2}\begin{pmatrix}0&X\\\ltrans{\ol{X}}&0\end{pmatrix}+\frac{1}{2}\begin{pmatrix}0&-iX\\ i\ltrans{\ol{X}}&0\end{pmatrix}\otimes \bfi, &X\in M_{n,n}(\bC),
\end{align*}   
to which we can apply the conjugation by $\fs$ and obtain the weight raising operators for $\U(J_{2n})$,
\begin{align*}
   \mu^+_{\U(J_{2n}),X}&=\fs^{-1} \mu^+_{\U(n,n),X}\,\fs.
\end{align*}
Denote by $\fU^+_{\U(n,n)}$ (resp. $\fU^+_{\U(J_{2n})}$) the $\bC$-vector space spanned by $\mu^+_{\U(n,n),X}$ (resp. $ \mu^+_{\U(J_{2n}),X}$), $X\in M_{n,n}(\bC)$.

We set $\bm{\mu}^{+,\upl}_{\U(n,n)}$ and $\bm{\mu}^{+,\upl}_{\U(J_{2n})}$ (resp. $\bm{\mu}^{+,\lowr}_{\U(n,n)}$ and $\bm{\mu}^{+,\lowr}_{\U(J_{2n})}$) to be the $a\times a$ (resp. $b\times b$) matrices taking values in $\fU^+_{\U(n,n)}$ (resp. $\fU^+_{\U(J_{2n})}$)  with the $(i,j)$-entry as
\begin{align*}
   \bm{\mu}^{+,\upl}_{\U(n,n),ij}&=\mu^+_{\U(n,n),E_{ij}}, &\bm{\mu}^{+,\upl}_{\U(J_{2n}),ij}&=\mu^+_{\U(J_{2n}),E_{ij}}\\
   \bm{\mu}^{+,\lowr}_{\U(n,n),ij}&=\mu^+_{\U(n,n),E_{a+i,a+j}}, &\bm{\mu}^{+,\lowr}_{\U(J_{2n}),ij}&=\mu^+_{\U(J_{2n}),E_{a+i,a+j}},
\end{align*}
where $E_{ij}$ denotes the $n\times n$ matrix with $1$ as the $(i,j)$-entry and $0$ elsewhere.

 The vector space spanned by the entries of $\bm{\mu}^{+,\upl}_{\U(n,n)}$ and $\bm{\mu}^{+,\lowr}_{\U(n,n)}$ is isomorphic to the quotient 
\[
   \left.\fU^+_{\U(n,n)}\right/\fU^+_{\U(n,n)}\cap \imath\big(\left(\Lie \U(a,b)\times\U(b,a)\right)\otimes_{\bR}\bC\big),
\] 
where $\imath$ is the embedding 
\begin{equation}\label{eq:embedding}
\begin{aligned}
   \imath:\U(a,b)\times\U(b,a)&\lhra \U(n,n)\\
   \begin{blockarray}{ccc}a&b\\\begin{block}{(cc)c}x_1&x_2&a\\x_3&x_4&b\\\end{block}\end{blockarray}\times \begin{blockarray}{ccc}b&a\\\begin{block}{(cc)c}y_1&y_2&b\\y_3&y_4&a\\\end{block}\end{blockarray}&\longmapsto \begin{blockarray}{ccccc}a&b&a&b\\\begin{block}{(cccc)c}x_1&&&x_2&a\\&y_1&y_2&&b\\&y_3&y_4&&a\\x_3&&&x_4&b\\\end{block}\end{blockarray}.
\end{aligned}
\end{equation}
The idea of considering this space for constructing archimedean sections for arithmetic applications of the doubling method stems from \cite{HaBun}.

For $k$, $r$, and $\left(\utau;\unu\right)$ such that both $\tau_a-\frac{k+r}{2}$ and $\nu^*_b-\frac{k-r}{2}$ are non-negative integers, define the polynomial $\fQ_{k, r,(\utau;\unu)}$ on $n\times n$ matrices by
\begin{equation}\label{eq:Deltaj}
   \fQ_{k,r,(\utau;\unu)}=\prod_{j=1}^{a-1}\Delta_j^{\tau_j-\tau_{j+1}}\Delta_a^{\tau_a-\frac{k+r}{2}}\prod_{j=1}^{b-1}{\Delta'_j}^{\nu^*_j-\nu^*_{j+1}}{\Delta'_b}^{\nu^*_b-\frac{k-r}{2}},
\end{equation}
where $\Delta_j$ (resp. $\Delta'_j$) stands for the determinant of the upper left (lower right) $j\times j$ block of a matrix, and $\unu^*=(\nu^*_1,\dots,\nu^*_b)=(-\nu_b,\dots,-\nu_1)$. (N.B. These are similar to the polynomials arising in the differential operators in \cite[Corollary 5.2.10]{EFMV} and \cite[Section 12]{shar}.) Note that the condition that $\tau_a-\frac{k+r}{2}$ and $\nu^*_b-\frac{k-r}{2}$ are non-negative integers is essentially the later defined condition \eqref{kcondition}, which is the condition for the point $s=\pm\frac{k-n}{2}+\frac{1}{2}$ to be critical for $L\left(s,\cD^*_{(\utau;\unu)}\times\chi^r_{\mr{ac}}\right)$.

We define the differential operators $D^{\U(n,n)}_{k,r,(\utau;\unu)}$ and $D^{\U(J_{2n})}_{k,r,(\utau;\unu)}$ by
\begin{equation}\label{eq:Ddef}
\begin{aligned}
   D^{\U(n,n)}_{k,r,(\utau;\unu)}&=\fQ_{k, r,(\utau;\unu)}\left(\frac{\bm{\mu}^{+,\lowr}_{\U(n,n)}}{2\pi i}\right), &D^{\U(J_{2n})}_{k,r,(\utau;\unu)}&=\fQ_{k, r,(\utau;\unu)}\left(\frac{\bm{\mu}^{+,\lowr}_{\U(J_{2n})}}{2\pi i}\right).
\end{aligned}
\end{equation}

The section $f_{k,(\utau;\unu)}(s,\chi^r_{\mr{ac}})\in I(s,\chi^r_{\mr{ac}})$ used in \cite{HELS} for constructing $p$-adic $L$-functions is defined by
\begin{equation}\label{eq:fktv}
\begin{aligned}
   f_{k,(\utau;\unu)}(s,\chi^r_{\mr{ac}}):=D^{\U(J_{2n})}_{k,r,(\utau;\unu)} f_k(s,\chi^r_{\mr{ac}}).
\end{aligned}
\end{equation}

\subsection{The archimedean doubling zeta integrals}
For a section $f(s,\chi^r_{\mr{ac}})\in I(s,\chi^r_{\mr{ac}})$, and vectors $v^*_1\in\cD^*_{(\utau;\unu)}$, $v_2\in \cD_{(\utau;\unu)}$, the archimedean doubling zeta integral is defined by
\begin{equation}\label{eq:zeta}
\begin{aligned}
   Z\left(f(s,\chi^r_{\mr{ac}}),v^*_1,v_2\right)=\int_{\U(a,b)}f(s,\chi^r_{\mr{ac}})\left(\fs^{-1}\imath(g,1) \fs\right)\,\left<g\cdot v^*_1,v_2\right>\,dg,
\end{aligned}
\end{equation}
where  $\fs$ is the matrix defined in Equation \eqref{eq:s} and $\imath$ is the embedding \eqref{eq:embedding}.

In the following, we compute
\begin{equation}\label{ZItoCompute}
   \left.Z\left(f_{k,(\utau;\unu)}(s,\chi^r_{\mr{ac}}),v^*_{(\utau;\unu)},v_{(\utau;\unu)}\right)\right|_{s=\pm\frac{k-n}{2}}
\end{equation}
for integers $k$ satisfying the conditions
\begin{equation}\label{kcondition}
   \begin{aligned}
   k&\equiv r\mod 2, &k&\geq n, &\tau_a\geq\frac{k+r}{2},\,\nu_1\leq -\frac{k-r}{2},
\end{aligned}
\end{equation}
and the section $f_{k,(\utau;\unu)}(s,\chi^r_{\mr{ac}})$ defined in Equation \eqref{eq:fktv}, $v_{(\utau;\unu)}\in\cD_{(\utau;\unu)}$ the highest weight vector in the lowest $K$-type of $\cD_{(\utau;\unu)}$, and $v^*_{(\utau;\unu)}\in\cD^*_{(\utau;\unu)}$ the dual vector of $v_{(\utau;\unu)}$.

\begin{rmk}
The term $Z\left(f_{k,(\utau;\unu)}\left(\frac{k-n}{2},\chi^r_{\mr{ac}}\right),v^*_{(\utau;\unu)},v_{(\utau;\unu)}\right)$ is left as an unknown number in \cite[Main Theorem 9.2.2]{HELS} (which notes it is equal to a nonzero rational number multiplied by an automorphic period).  This paper computes this number, and the result is stated in Theorem~\ref{thm:Zright}.  
\end{rmk}

\section{Theta correspondence of unitary groups}\label{theta-section}
Our computation of the archimedean doubling zeta integral \eqref{ZItoCompute} relies on results from the theory of the theta correspondence between $\U(a,b)$, $\U(n,n)$ and the compact unitary group $\U(k)$. Before starting the computation, we briefly recall some basics about Weil representations.

\subsection{The Schr\"{o}dinger model}\label{section:schrodinger}
We fix the additive character 
\begin{align*}
   \be_{\bR}:\bR&\lra\bC^\times, &\be_{\bR}(x)=e^{2\pi i\cdot x}.
\end{align*}
For a positive integer $m$, let $\Sp(2m)$ be the symplectic group
$$
   \left\{g\in\GL(2m,\bR)\middle|\,\ltrans{g}\begin{pmatrix}&\bid_m\\-\bid_m&\end{pmatrix}g=\begin{pmatrix}&\bid_m\\-\bid_m&\end{pmatrix}\right\},
$$
and let $\Mp(2m)$ be the central extension
$$
   1\lra \bC^\times\lra \Mp(2m)\lra \Sp(2m)\lra 1.
$$
The metaplectic group $\Mp(2m)$ has a Weil representation. The definition of $\Mp(2m)$ and its Weil representation depend on the choice of the additive character of $\bR$.

Let $\sS_{mk}$ denote the Schr\"{o}dinger model of the Weil representation of $\Mp(2mk)$.  (See \cite{kudla-notes}, especially Chapters I and II for an introduction to the Schr\"odinger model.)  More precisely, $\sS_{mk}$ is the space of Schwartz functions on $\bR^{mk}$ (valued in $\bC$) on which $\left(g=\begin{pmatrix}A&B\\C&D\end{pmatrix},z\right)\in\Mp(2mk)$ acts by sending the Schwartz function $\phi(X)\in\sS_{mk}$ to
\begin{equation}\label{eq:weiloper}
\begin{aligned}
   X\longmapsto &\chi^{\be_{\bR}}_k(x(g),z)\,\gamma_{\bR}(\be_{\bR})^{-\mr{rk} C}\\
   &\times\int_{\bR^{mk}/\ker \ltrans{C}}\phi\left(\ltrans{C}Y+\ltrans{A}X\right)\,\be_\bR\left(\frac{1}{2}\left(\ltrans{Y}C\ltrans{D}Y+\ltrans{Y}C\ltrans{B}X+\ltrans{X}A\ltrans{B}X\right)\right)\,d_gY,
\end{aligned}
\end{equation}
where $\chi^{\be_{\bR}}_k$ is the character 
\begin{align*}
    \bR^\times\ltimes\bC^\times&\lra\bC^\times\\
    (x,z)&\longmapsto \left(x,(-1)^{\frac{mk(mk-1)}{2}}\right)_{\bR}\cdot\begin{cases}1,&\text{if $k$ is even},\\ z\cdot\gamma_{\bR}(x,\be_{\bR})^{-1},&\text{if $k$ is odd},\end{cases}
\end{align*}
with $(, )_\bR$ denoting the Hilbert symbol, $\gamma_{\bR}(\be_{\bR})=e^{\frac{\pi i}{4}}$, $\gamma_{\bR}(x,\be_{\bR})=e^{\frac{\pi i}{4}(\sgn (x)-1)}$,  $x(g)\in\bR$ is defined by
\begin{align*}
   x(g)&=\det(A_1A_2), &g&=\begin{psm}A_1&B_1\\&\ltrans{A}^{-1}_1\end{psm}\begin{psm}\bid_{mk-\mr{rk}C}&&&\\&&&\bid_{\mr{rk}C}\\&&\bid_{mk-\mr{rk}C}&\\&-\bid_{\mr{rk}C}&\end{psm}\begin{psm}A_2&B_2\\&\ltrans{A}^{-1}_2\end{psm},
\end{align*}
and the Haar measure $d_gY$ on $\bR^{mk}/\ker C$ is defined such that operator \eqref{eq:weiloper} preserves the $L^2$-norm on $\sS_{mk}$.

\subsection{The Siegel--Weil sections}
We consider the reductive dual pair $(\U(J_{2n}),\U(k))$ inside $\Sp(4nk)$. Thanks to our assumption that $k$ and $r$ have the same parity, we can fix a lifting
\begin{equation}\label{eq:iotar}
\xymatrix{
    \U(J_{2n})\times\U(k)\,\ar@{^{(}->}[r]^-{\iota_r}\ar@{->}[dr]&\Mp(4nk)\ar[d]\\
    &\Sp(4nk)
}
\end{equation}
and an isomorphism of $\bR^{2nk}$ with $M_{n,k}(\bC)$, such that the action of $\U(J_{2n})\times\U(k)$ on $\sS_{2nk}$ satisfies  
\begin{align*}
   \omega_{\sS_{2nk}}\left(\iota_r\left(\begin{psm}A&B\\&\ltrans{\bar{A}}^{-1}\end{psm},1\right)\right)\Phi(X)&= \chi^r_{\mr{ac}}(\det A)|\det A|^{\frac{k}{2}}_{\bC}\,\be_{\bR}\left(\frac{1}{2} \Tr\ltrans{\bar{X}}B\ltrans{\bar{A}}X\right)\Phi(\ltrans{\bar{A}}X), 
\end{align*}
where $X\in M_{n,k}(\bC)$ and $\Phi\in\sS_{2nk}$ is viewed as a Schwartz function on $M_{n,k}(\bC)\cong \bR^{2nk}$. In the following, for $(g,h)\in \U(J_{2n})\times\U(k)$ and $\iota_r$ as in the diagram \eqref{eq:iotar}, we denote $\omega_{\sS_{2nk}}\left(\iota_r(g,h)\right)$ by $\omega_{\sS_{2nk}}(g,h)$.

The Siegel--Weil section $f_{\mr{SW}}(\Phi)\in  I\left(\frac{k-n}{2},\chi^r_{\mr{ac}}\right)$ associated to $\Phi\in\sS_{2nk}$ is defined by
$$
   f_{\mr{SW}}(\Phi)(g):=\omega_{\sS_{2nk}}(g,1)\Phi(0).
$$
Let $\Phi_0\in\sS_{2nk}$ be the Gaussian function defined by
$$
\Phi_0(X):=e^{-\pi\Tr \ltrans{\bar{X}}X}.
$$
Then it follows easily from the formula \eqref{eq:weiloper} that the evaluation at $s=\frac{k-n}{2}$ of the classical section
 $f_k(s,\chi^r_{\mr{ac}})$ (given in Equation \eqref{eq:fk}) is the Siegel--Weil section attached to the Gaussian function, i.e.
\begin{equation}\label{eq:SWGaussian}
   f_{\mr{SW}}(\Phi_0)=f_k\left(\frac{k-n}{2},\chi^r_{\mr{ac}}\right).
\end{equation}

\subsection{Restrictions of Siegel--Weil sections and matrix coefficients of Weil representations}
Let $V$ (resp. $V'$) be an $n$-dimensional Hermitian space over $\bC$ with signature $(a,b)$ (resp. $(b,a)$), and set $\bV=V\oplus V'$. We fix the following basis
\begin{align}
   V&:(\ue,\uf)=(e_1,\dots,e_a,f_1,\dots,f_b)\label{eq:bV}\\
   V'&:(\ue',\uf')=(e'_1,\dots,e'_b,f'_1,\dots,f'_a)\label{eq:bV'}\\
   \bV&:(\ue,\ue',\uf',\uf),\label{eq:bV3}
\end{align}
under which the matrix of the Hermitian form on $V$ (resp. $V'$) is $\begin{pmatrix}\bid_a&\\&-\bid_b\end{pmatrix}$ (resp. $\begin{pmatrix}\bid_b&\\&-\bid_a\end{pmatrix}$), and the embedding $\U(V)\times\U(V')\hra \U(\bV)$ agrees with the embedding $\imath$ in \eqref{eq:embedding}.

We consider another basis of $\bV$ given as
\begin{equation}\label{eq:bV4}
   (\ue,\ue',\uf',\uf)\,\fs=\left(\frac{\ue+\uf'}{\sqrt{2}},\,\frac{\ue'+\uf}{\sqrt{2}},\,\frac{-i(\ue-\uf')}{\sqrt{2}},\,\frac{-i(\ue'-\uf)}{\sqrt{2}}\right).
\end{equation}
With respect to the basis \eqref{eq:bV3}, $\U(\bV)$ is identified with $\U(n,n)$, and with respect to \eqref{eq:bV4}, $\U(\bV)$ is identified with $\U(J_{2n})$.

We also consider the symplectic spaces over $\bR$ associated to the Hermitian spaces over $\bC$. Denote by $V_{\bR}$ (resp. $V'_{\bR}$, $\bV_{\bR}$) the symplectic space over $\bR$ of dimension $2n$ (resp. $2n$, $4n$), which equals $V$ (resp. $V'$, $\bV$) viewed as an $\bR$-vector space and equipped with the symplectic form obtained by taking the imaginary part of the Hermitian form. For $V_{\bR}$ (resp. $V'_{\bR}$), we fix the following basis
\begin{align}
   V_{\bR}&:(\ue,\uf,-i\ue,i\uf),\\
   V'_{\bR}&:(\uf',\ue',i\uf',-i\ue'),\\
   \bV_{\bR}&:(\ue,\ue',\uf',\uf,-i\ue,-i\ue',i\uf',i\uf),\label{eq:bV1}\\
   \bV_{\bR}&:(\ue,\ue',\uf',\uf,-i\ue,-i\ue',i\uf',i\uf)\,\fS,\label{eq:bV2}
\end{align}
where $\fS=\frac{1}{\sqrt{2}}\begin{pmatrix}\bid_n&0&0&\bid_n\\\bid_n&0&0&-\bid_n\\0&-\bid_n&\bid_n&0\\0&\bid_n&\bid_n&0\end{pmatrix}$. The bases \eqref{eq:bV4} and \eqref{eq:bV2} are compatible with our fixed embedding \eqref{eq:iotar}.

With these fixed bases, we obtain the following commutative diagram:
$$\xymatrixcolsep{5pc}\xymatrix{
   \U(a,b)\times\U(b,a)\,\ar@{^{(}->}[r]^-{\imath}\ar@{^{(}->}[d]_{\iota^\natural_{r,1}\times\iota^\natural_{r,2}}&\U(n,n)\ar[r]^-{\fs^{-1}\cdot\,\fs}\ar@{^{(}->}[d]^{\iota^\natural_r}&\U(J_{2n})\ar@{^{(}->}[d]^{\iota_r}\\
   \Mp(2n)\times\Mp(2n)\,\ar@{^{(}->}[r]&\Mp(4n)\ar[r]^-{\fS^{-1}\cdot\,\fS}&\Mp(4n)
}.
$$
The action of $\U(a,b)$ (resp. $\U(b,a)$) on $\sS_{nk}$ is induced by the embedding $\iota^\natural_{r,1}$ (resp. $\iota^\natural_{r,2}$) and the standard action of $\Mp(2n)$ on $\sS_{nk}$. Similarly, we make  $\U(n,n)$ act on $\sS_{2nk}$ through the embedding $\iota^\natural_r$ and the standard action of $\Mp(4n)$ on $\sS_{2nk}$. 

The embedding $\iota^\natural_{r,1}$ (resp. $\iota^\natural_{r,2}$) also induces embedding of the reductive dual pair $\U(a,b)\times\U(k)$ (resp.$\U(a,b)\times\U(k)$) into $\Mp(2nk)$, and we denote by $\omega_{\sS_{nk}}(g,h)$ the action of $(g,h)\in \U(a,b)\times\U(k)$ on $\sS_{nk}$. We define the matrix coefficients for the action of $\U(a,b)$ on $\sS_{nk}$ as
\begin{equation}\label{eq:SMC}
\begin{aligned}
   \mr{MC}_{\sS_{nk}}:\U(a,b)\times (\sS_{nk}\otimes\sS_{nk})&\lra\bC\\
   (g,\phi_1\otimes\phi_2)&\longmapsto \int_{\bR^{nk}}\omega_{\sS_{nk}}(g,1)\phi_1(x)\,\phi_2(x)\,dx.
\end{aligned}
\end{equation}
Note that $\sS_{nk}\otimes\sS_{nk}$ is dense in $\sS_{2nk}$. By taking the limit, one can extend the above map to 
$$\mr{MC}_{\sS_{nk}}:\U(a,b)\times\sS_{2nk}\lra\bC.$$ 

Let $\wt{\fS}=[\fS,1]\in \Mp(4n)$ and $z_{\fS}\in \bC^\times$ such that $\wt{\fS}^{-1}=[\fS^{-1},z_{\fS}]$.
\begin{prop}\label{prop:MC}
For $\phi_1\otimes\phi_2\in\sS_{nk}\otimes\sS_{nk}$, 
$$
   \omega_{\sS_{2nk}}\left(\wt{\fS}^{-1},1\right)(\phi_1\otimes\phi_2)(0)=\chi^{\be_{\bR}}_k(1,z_{\fS})\,\gamma_{\bR}(\be_{\bR})^{-nk}\, 2^{\frac{nk}{2}}\int_{\bR^{nk}}\phi_1(x)\,\phi_2(-x)\,dx.
$$
\end{prop}
\begin{proof}
We have 
\[\fS^{-1}=\frac{1}{\sqrt{2}}\begin{pmatrix}\bid_n&\bid_n&0&0\\0&0&-\bid_n&\bid_n\\0&0&\bid_n&\bid_n\\\bid_n&-\bid_n&0&0\end{pmatrix}=\protect\begin{pmatrix}\bid_n&0&0&0\\0&0&0&\bid_n\\0&0&\bid_n&0\\0&-\bid_n&0&0\protect\end{pmatrix}\frac{1}{\sqrt{2}}\protect\begin{pmatrix}\bid_n&\bid_n&0&0\\-\bid_n&\bid_n&0&0\\0&0&\bid_n&\bid_n\\0&0&-\bid_n&\bid_n\protect\end{pmatrix}.
\] 
It follows from the formula \eqref{eq:weiloper} that $x(\fS^{-1})=1$, and
\begin{align*}
   &\omega_{\sS_{2nk},1}\left(\wt{\fS}^{-1},1\right)(\phi_1\otimes\phi_2)(0)\\
   =&\,\chi^{\be_\bR}_k(1,z_{\fS})\,\gamma_{\bR}(\be_\bR)^{-nk}\\
   &\times\int_{M_{n,k}(\bR)}(\phi_1\otimes\phi_2)\left(\frac{1}{\sqrt{2}}\begin{pmatrix}0&\bid_n\\0&-\bid_n\end{pmatrix}\begin{pmatrix}0\\x\end{pmatrix}\right)\,\be_{\bR}\left(\frac{1}{4}\,\Tr \begin{pmatrix}0&\bid_n\\0&-\bid_n\end{pmatrix}\begin{pmatrix}0\\x\end{pmatrix}\begin{pmatrix}0&\ltrans{x}\end{pmatrix}\begin{pmatrix}\bid_n&\bid_n\\0&0\end{pmatrix}\right)\,dx\\
   =&\,\chi^{\be_\bR}_k(1,z_{\fS})\,\gamma_{\bR}(\be_\bR)^{-nk}\, 2^{\frac{nk}{2}}\int_{M_{n,k}(\bR)}\phi_1(x)\,\phi_2(-x)\,dx.
\end{align*}
\end{proof}

The next proposition relates the evaluation at $s=\frac{k-n}{2}$ of the section defined in \eqref{eq:fktv} to the matrix coefficient of the Schwartz function obtained by applying the differential operator defined in \eqref{eq:Ddef} to the Gaussian function. 
\begin{prop}\label{prop:f-MC}
Maintaining the conventions from above, we have
\begin{align*}
   &f_{k,(\utau;\unu)}\left(\frac{k-n}{2},\chi^r_{\mr{ac}}\right)\left(\fs^{-1}\imath(g,1)\fs\right)\\
   =&\,(-1)^{\frac{a(k+r)}{2}+\frac{b(k-r)}{2}+\sum\tau_j+\sum\nu_j}2^{\frac{nk}{2}}\,\mr{MC}_{\sS_{nk}}\left(g,\omega_{\sS_{2nk}}\left(D^{\U(n,n)}_{k,r,(\utau;\unu)}\right)\,\Phi_0\right).
\end{align*}
\end{prop}
\begin{proof}
It follows from the definition of $f_{k,(\utau;\unu)}\left(\frac{k-n}{2},\chi^r_{\mr{ac}}\right)$ and \eqref{eq:SWGaussian} that
\begin{equation}\label{eq:fMC}
\begin{aligned}
   &f_{k,(\utau;\unu)}\left(\frac{k-n}{2},\chi^r_{\mr{ac}}\right)\left(\fs^{-1}\imath(g,1)\fs\right)\\
   =&\,\omega_{\sS_{2nk}}\left(\wt{\fS}^{-1},1\right)\omega_{\sS_{2nk}}\left(\imath(g,1)\right)\omega_{\sS_{2nk}}\left(\wt{\fS},1\right)\omega_{\sS_{2nk}}\left(D^{\U(J_{2n})}_{k,r,(\utau;\unu)}\right)\Phi_0(0)
\end{aligned}
\end{equation}
By the formulas \eqref{eq:weiloper} and
\[ 
  \fS=\frac{1}{\sqrt{2}}\protect\begin{pmatrix}\bid_n&0&0&\bid_n\\\bid_n&0&0&-\bid_n\\0&-\bid_n&\bid_n&0\\0&\bid_n&\bid_n&0\protect\end{pmatrix}=\frac{1}{\sqrt{2}}\protect\begin{pmatrix}\bid_n&\bid_n&0&0\\\bid_n&-\bid_n&0&0\\0&0&\bid_n&\bid_n\\0&0&\bid_n&-\bid_n\protect\end{pmatrix}\protect\begin{pmatrix}\bid_n&0&0&0\\0&0&0&\bid_n\\0&0&\bid_n&0\\0&-\bid_n&0&0\protect\end{pmatrix},
\] 
we see that $x(\fS)=(-1)^n$ and 
\begin{align*}
   &\omega_{\sS_{2nk}}(\wt{\fS},1)\,\Phi_0\begin{pmatrix}x_1\\x_2\end{pmatrix}\\
   =&\,\chi^{\be_\bR}_k((-1)^n,1)\gamma_{\bR}(\be_\bR)^{-nk}\,2^{-\frac{nk}{2}}\int_{M_{n,k}(\bR)}\Phi_0\begin{pmatrix}\frac{x_1+x_2}{\sqrt{2}}\\-\sqrt{2}y\end{pmatrix}\be_{\bR}\left(\Tr\ltrans{y}(x_2-x_1)\right)\,dy\\
   =&\,\chi^{\be_\bR}_k((-1)^n,1)\gamma_{\bR}(\be_{\bR})^{-nk}\Phi_0\begin{pmatrix}x_1\\x_2\end{pmatrix}.
\end{align*} 
Combining this with Proposition~\ref{prop:MC}, we get
\begin{equation}\label{eq:MC1}
\begin{aligned}
   \eqref{eq:fMC}=&\,\omega_{\sS_{2nk}}\big(\wt{\fS}^{-1},1\big)\omega_{\sS_{2nk}}\left(\imath(g,1)\right)\omega_{\sS_{2nk}}\left(D^{\U(n,n)}_{k,r,(\utau;\unu)}\right)\omega_{\sS_{2nk}}\big(\wt{\fS},1\big)\Phi_0(0)\\
   =&\,\chi^{\be_\bR}_k((-1)^n,1)\chi^{\be_\bR}_k(1,z_{\fS})\gamma_\bR(\be_\bR)^{-2nk}\,2^{\frac{nk}{2}}\\
   &\times\int_{\bR^{nk}}\omega_{\sS_{nk}}(g,\bid_k)\omega_{\sS_{2nk}}\left(D^{\U(n,n)}_{k,r,(\utau;\unu)}\right)\Phi_0\begin{pmatrix}x\\-x\end{pmatrix}\,dx.
\end{aligned}
\end{equation}
In particular, this equality holds when the differential operator $D^{\U(J_{2n})}_{k,r,(\utau;\unu)}$ is the trivial one and $g=\bid_n$. This case implies that $\chi^{\be_\bR}_k((-1)^n,1)\chi^{\be_\bR}_k(1,z_{\fS})\gamma_\bR(\be_\bR)^{-2nk}=1$. Therefore,
\begin{equation}\label{eq:MC2}
   \eqref{eq:fMC}=2^{\frac{nk}{2}}\int_{\bR^{nk}}\omega_{\sS_{nk}}(g,1)\omega_{\sS_{2nk}}\left(D^{\U(n,n)}_{k,r,(\utau;\unu)}\right)\Phi_0\begin{pmatrix}x\\-x\end{pmatrix}\,dx.
\end{equation}

For the action of $\U(b)\times\U(a)\subset\U(b,a)$ through the Weil representation, the Gaussian function $\Phi_0$ is of highest weight $\left(\frac{k-r}{2},\dots,\frac{k-r}{2};\frac{k+r}{2},\dots,\frac{k+r}{2}\right)$. Then it follows from the definition of the differential operator $D^{\U(n,n)}_{k,r,(\utau,\unu)}$ in \eqref{eq:Ddef} that $D^{\U(n,n)}_{k,r,(\utau,\unu)}\Phi_0$ is of highest weight $(\utau;\unu)$. Thus,
\begin{equation}\label{eq:MC3}
   \omega_{\sS_{2nk}}\left(D^{\U(n,n)}_{k,r,(\utau,\unu)}\right)\Phi_0(x_1,-x_2)=(-1)^{\frac{a(k+r)}{2}+\frac{b(k-r)}{2}+\sum\tau_j+\sum\nu_j}\omega_{\sS_{2nk}}\left(D^{\U(n,n)}_{k,r,(\utau,\unu)}\Phi_0\right)(x_1,x_2).
\end{equation}
Combining Equations \eqref{eq:MC2} and \eqref{eq:MC3}, we get
\begin{align*}
   \eqref{eq:fMC}=(-1)^{\frac{a(k+r)}{2}+\frac{b(k-r)}{2}+\sum\tau_j+\sum\nu_j}\,2^{\frac{nk}{2}}\,\mr{MC}_{\sS_{nk}}\left(g,\omega_{\sS_{2nk}}\left(D^{\U(n,n)}_{k,r,(\utau;\unu)}\right)\,\Phi_0\right).
\end{align*}

\end{proof}

\begin{rmk}
Viewing the evaluations at $\fs^{-1}\imath(g,1)\fs$ of the Siegel--Weil sections as the matrix coefficients of Weil representations has been widely used in the study of the theta correspondence and the doubling method, for example in \cite{LiThetaCoh}. We write out all the computation above in order to make sure that the factors in the comparison of the evaluations of the sections and the matrix coefficients are precise.
\end{rmk}

\subsection{The Fock model}\label{section:fock}
Proposition~\ref{prop:f-MC} allows us to reduce the computation of \eqref{ZItoCompute} in the case $s=\frac{k-n}{2}$ to studying the decomposition of the Weil representation of $\U(a,b)\times\U(k)$ and its matrix coefficients. Unlike the case of symplectic groups considered in \cite{liu-archimedean}, the Schr{\"o}dinger model is not convenient for this purpose due to the lack of nice explicit formulas for the action of $\U(a,b)$ and pluri-harmonic polynomials (especially when $a\neq b$). (Note that the model used in \cite{KVWeil} is not the Schr{\"o}dinger model when $a\neq b$.) We need to introduce the Fock model and the Bargmann transform. 

Following \cite{Folland}, let $\sF_{mk}$ be the Fock model of the Weil representation of the metaplectic group $\Mp(2mk)$. It consists of entire functions on $\bC^{mk}$ which are square integrable with respect to the the Hermitian pairing
\begin{equation}\label{eq:Fpairing}
   \left<F_1,F_2\right>_{\sF_{mk}}=2^{-mk}\int_{\bC^{mk}}F_1(z)\ol{F_2(z)}\,e^{-\pi\ltrans{\ol{z}}z}\,|dz\,d\ol{z}|.
\end{equation}
Let
$$
   \cW=\frac{1}{\sqrt{2}}\begin{pmatrix}\bid_{nk}&\bfi\cdot\bid_{nk}\\-\bid_{nk}&\bfi\cdot\bid_{nk}\end{pmatrix}.
$$
For $g\in\Sp(2mk)$ with $\cW g\cW^{-1}=\begin{pmatrix}P&Q\\\ol{Q}&\ol{P}\end{pmatrix}$, its action on $F\in\sF_{mk}$ is given by
\begin{equation}\label{eq:Fock}
\begin{aligned}
   &\omega_{\sF_{mk}}(g)F(Z)\\
   =&\,\left(\mr{det}^{-1/2} P\right)\, 2^{-mk}\int_{\bC^{mk}}e^{\frac{\pi}{2}(\ltrans{Z}\ol{Q}P^{-1}Z+2\ltrans{\ol{W}}P^{-1}\ltrans{Z}-\ltrans{\ol{W}}P^{-1}Q\ol{W})}\,F(W)e^{-\pi \ltrans{W}\ol{W}}\,|dWd\ol{W}|.
\end{aligned}
\end{equation}
Here the ambiguity of $\mr{det}^{1/2} P$ is because the Weil representation is a representation of $\Mp(2mk)$ rather than $\Sp(2mk)$, but for our purpose it is not necessary to be precise about this ambiguity. The Hermitian pairing \eqref{eq:Fpairing} is equivaraint for this action.  From \eqref{eq:Fock}, it is easy to deduce formulas for the action of $\left(\Lie\Sp(2mk)\right)\otimes_{\bR}\bC$ on $\sF_{2mk}$. For $X\in\Sym(mk,\bR)$, let
\begin{align*}
   \mu^+_{\Sp(2mk),X}&=\frac{1}{2}\begin{pmatrix}0&X\\X&0\end{pmatrix}+\frac{1}{2}\begin{pmatrix}-X&0\\0&X\end{pmatrix}\otimes i, &\mu^-_{\Sp(2mk),X}&=\frac{1}{2}\begin{pmatrix}0&X\\X&0\end{pmatrix}+\frac{1}{2}\begin{pmatrix}X&0\\0&-X\end{pmatrix}\otimes i.
\end{align*}
Since
\begin{align*}
   \cW\mu^+_{\Sp(2mk),X}\cW^{-1}&=\begin{pmatrix}0&0\\X&0\end{pmatrix}\otimes i, &\cW\mu^+_{\Sp(2mk),X}\cW^{-1}&=-\begin{pmatrix}0&X\\0&0\end{pmatrix}\otimes i,
\end{align*}
by the formulas \eqref{eq:Fock}, we easily see that
\begin{align}
   \label{eq:Fwtr}\omega_{\sF_{mk}}\left(\mu^+_{\Sp(2mk),X}\right)F(Z)&=\frac{\pi i}{2}\,\ltrans{Z}XZ\cdot F(Z), \\
   \label{eq:Fwtl}\omega_{\sF_{mk}}\left(\mu^-_{\Sp(2mk),X}\right)F(Z)&=\frac{i}{2\pi}\sum_{1\leq i,j\leq mk} X_{ij}\frac{\partial^2}{\partial Z_{i}\partial Z_{j}}F(Z).
\end{align}

Between the Sch\"{o}dinger model and the Fock model, there is the Bargmann transform
\begin{equation}\label{eq:Bargmann}
\begin{aligned}
   \sB:\sS_{mk}&\lra \sF_{mk}\\
   \phi&\longmapsto \sB(\phi)(Z)=2^{\frac{mk}{4}}\int_{\bR^{nk}}f(X)\,e^{2\pi X\ltrans{Z}-\pi X\ltrans{X}-\frac{\pi}{2}Z\ltrans{Z}}\,dX.
\end{aligned}
\end{equation}
The Bargmann transform is an isometry with respect to the standard Hermitian pairing on $\sS_{mk}$ and the Hermitian pairing \eqref{eq:Fpairing},  i.e.
$$
   \int_{\bR^{mk}}\phi_1(X)\,\ol{\phi_2(X)}\,dX=2^{-mk}\int_{\bC^{mk}}\sB(\phi_1)(Z)\ol{\sB(\phi_2)(Z)}\,e^{-\pi\ltrans{\ol{Z}}Z}\,|dZd\ol{Z}|.
$$
Similarly to the definition of the matrix coefficients in \eqref{eq:SMC}, we define the matrix coefficients for the action of $\U(a,b)$ on $\sF_{nk}$ as
\begin{align*}
   \mr{MC}_{\sF_{nk}}:\U(a,b)\times (\sF_{nk}\otimes\sF_{nk})&\lra\bC\\
   (g,F_1\otimes F_2)&\longmapsto 2^{-nk}\int_{\bC^{nk}}\omega_{\sF_{nk}}(g,1)F_1(Z)\,F_2(\ol{Z})\,e^{-\pi\ltrans{\ol{Z}}Z}\,|dZd\ol{Z}|,
\end{align*}
and extend it to $\U(a,b)\times\sF_{2nk}$. It follows from the isometry property of the Bargmann transform that
\begin{equation}\label{eq:SFMC}
   \mr{MC}_{\sS_{nk}}(g,\Phi)=\mr{MC}_{\sF_{nk}}(g,\sB(\Phi)).
\end{equation}

We identify $\bC^{nk}$ (resp. $\bC^{nk}$, $\bC^{2nk}$) with $M_{n,k}(\bC)$ (resp. $M_{n,k}(\bC)$, $M_{2n,k}(\bC)$) through the basis \eqref{eq:bV}\eqref{eq:bV3} and the standard basis of the positive definite Hermitain space of dimension $k$ over $\bC$ (for which the matrix of the Hermitian form is $\bid_k$). We write $z,w\in M_{n,k}(\bC)$ as
\begin{align}
   z&=\begin{blockarray}{cc}k\\\begin{block}{(c)c}z_1&a\\z_2&b\\\end{block}\end{blockarray},&w&=\begin{blockarray}{cc}k\\\begin{block}{(c)c}w_1&b\\w_2&a\\\end{block}\end{blockarray}.\label{eq:zcoor}
\end{align}
The space $M_{n,k}(\bC)\times M_{n,k}(\bC)$ is identified with $M_{2n,k}(\bC)$ by 
$$
   (z,w)\longmapsto \begin{pmatrix}z_1\\w_1\\w_2\\z_2\end{pmatrix}.
$$
Also, 
$$
   \cW\,\iota^\natural_r\left(\mu^+_{\U(n,n),X}\right)\cW^{-1}=-\begin{pmatrix}0&0&0&0\\0&0&0&0\\0&X&0&0\\\ltrans{X}&0&0&0\end{pmatrix}.
$$
Hence,
\begin{equation}\label{eq:Udiff}
   \omega_{\sF_{2nk}}\left(\mu^+_{\U(n,n),X}\right) F(z,w)=-\pi \, \Tr\left( X\begin{pmatrix}z_1\ltrans{w}_2&z_1\ltrans{z}_2\\w_1\ltrans{w}_2&w_1\ltrans{z}_2\end{pmatrix}\right) F(z,w).
\end{equation}
It is also easy to see that $h\in\U(k)$ acts on $F\in\sF_{2nk}$ by
\begin{equation}\label{eq:U(k)act}
    \omega_{\sF_{2nk}}(h)F(z,w)=F\left(\begin{pmatrix}z_1h\\z_2\ol{h}\end{pmatrix},\begin{pmatrix}w_1h\\w_2\ol{h}\end{pmatrix}\right).
\end{equation}

\begin{prop}\label{prop:F}
We have
\begin{align*}
   f_{k,(\utau;\unu)}\left(\frac{k-n}{2},\chi^r_{\mr{ac}}\right)\left(\fs^{-1}\imath(g,1)\fs\right)=&\,(2i)^{-\sum\tau_j-\sum\nu^*_j+\frac{a(k+r)}{2}+\frac{b(k-r)}{2}}\,\mr{MC}_{\sF_{nk}}\left(g,F_{k,r,(\utau;\unu)}\right),
\end{align*}
where
\begin{align*}
   &F_{k,r,(\utau;\unu)}(z,w)=\fQ_{k,r,(\utau;\unu)}\begin{pmatrix}z_1\ltrans{w}_2&z_1\ltrans{z}_2\\w_1\ltrans{w}_2&w_1\ltrans{z}_2\end{pmatrix}\\
   =&\prod_{j=1}^{a-1}\Delta_j\left(z_1\ltrans{w}_2\right)^{\tau_j-\tau_{j+1}}\det\left(z_1\ltrans{w}_2\right)^{\tau_a-\frac{k+r}{2}}\prod_{j=1}^{b-1}\Delta'_j\left(z_2\ltrans{w}_1\right)^{\nu^*_j-\nu^*_{j+1}}\det\left(z_2\ltrans{w}_1\right)^{\nu^*_b-\frac{k-r}{2}}.
\end{align*}
\end{prop}
\begin{proof}
Combining Proposition~\ref{prop:f-MC} and Equation \eqref{eq:SFMC}, we get
\begin{align*}
   f_{k,(\utau;\unu)}\left(\frac{k-n}{2},\chi^r_{\mr{ac}}\right)\left(\fs^{-1}\imath(g,1)\fs\right)= &\,(-1)^{\frac{a(k+r)}{2}+\frac{b(k-r)}{2}+\sum\tau_j+\sum\nu_j}2^{-\frac{nk}{2}}\\
   &\times\mr{MC}_{\sF_{nk}}\left(g,\omega_{\sF_{2nk}}\left(D^{\U(n,n)}_{k,r,(\utau;\unu)}\right)\sB(\Phi_0)\right).
\end{align*}
Direct computation shows that
\begin{equation*}
   \sB\left(\Phi_0\right)=2^{-\frac{nk}{2}}.
\end{equation*}
By Equation \eqref{eq:Udiff} and the definition of the operator $D^{\U(n,n)}_{k,r,(\utau;\unu)}$ in Equations \eqref{eq:Ddef}, we have
$$
   \omega_{\sF_{2nk}}\left(D^{\U(n,n)}_{k,r,(\utau;\unu)}\right) \sB(\Phi_0)=2^{-\frac{nk}{2}}\left(\frac{i}{2}\right)^{\sum\tau_j+\sum\nu^*_j-\frac{a(k+r)}{2}-\frac{b(k-r)}{2}}\fQ_{k,r,(\utau;\unu)}\left(\begin{pmatrix}z_1\ltrans{w}_2&z_1\ltrans{z}_2\\w_1\ltrans{w}_2&w_1\ltrans{z}_2\end{pmatrix}\right),
$$
and the proposition follows.
\end{proof}

Let $\sF^\circ_{nk}\subset\sF_{nk}$ be the subspace consisting of all the polynomials on $M_{n,k}(\bC)$. It is closed under the action of the maximal compact subgroup as well as the Lie algebra of $\Sp(2nk)$. By \cite[Theorem III(7.2)]{KVWeil}, the Weil representation of $\Lie\left(\U(a,b)\times\U(k)\right)$ on $\sF^\circ_{nk}$ ($k\geq n$) has a multiplicity free decomposition
\begin{equation}\label{eq:KVdecomp}
   \bigoplus_{\underline{m}\geq \frac{k+r}{2},\,\underline{l}^*\geq \frac{k-r}{2}}\cD_{(\underline{m};\underline{l})}\boxtimes \lambda_{k,r}(\um,\ul),
\end{equation}
where $\underline{m}$ (res. $\underline{l}$) is a dominant weight for $\GL(a)$ (resp. $\GL(b)$), and $\lambda_{k,r}(\um,\ul)$ denotes the irreducible representation of $\U(k)$ of highest weight
$$
   \left(m_1-\frac{k+r}{2},\dots,m_a-\frac{k+r}{2},0,\dots,0,l_1+\frac{k-r}{2},\dots,l_b+\frac{k-r}{2}\right).
$$
Denote by $\fH_{(a,b),k}$ the subspace of $\sF_{nk}$ consisting of joint harmonic polynomials for the action of $\U(a,b)$, i.e. polynomials annihilated by 
\begin{align*}
    \begin{pmatrix}0&X\\\ltrans{\ol{X}}&0\end{pmatrix}-\begin{pmatrix}0&-iX\\ i\ltrans{\ol{X}}&0\end{pmatrix}\otimes \bfi &\in\mu^-_{\U(a,b),X}, &X&\in M_{a,b}(\bC),
\end{align*}
or equivalently annihilated by
\begin{align*}
   &\sum_{j=1}^k\frac{\partial^2}{\partial z_{1,mj}\partial z_{2,lj}}, &\text{for all }1\leq m\leq a,\,1\leq l\leq b,
\end{align*}
where $z$ is the coordinate of $M_{n,k}(\bC)$ as in Equations \eqref{eq:zcoor}. (The equivalence follows easily from the formula in Equation \eqref{eq:Fwtl}.)  The joint harmonic polynomials correspond to vectors inside the direct sum of lowest $K$-types of the irreducible components in the decomposition \eqref{eq:KVdecomp}. The space $\fH_{(a,b),k}$ is preserved by the action of $\U(a)\times\U(b)\times\U(k)$, the maximal compact subgroup of $\U(a,b)\times\U(k)$, and has a decomposition 
$$\bigoplus\limits_{\um\geq\frac{k+r}{2},\,\ul^*\geq\frac{k-r}{2}}\fH_{(a,b),k}(\um,\ul)$$
with $\fH_{(a,b),k}(\um,\ul)$ corresponding to the lowest $K$-type of $\cD_{(\um;\ul)}\boxtimes\lambda_{k,r}(\um,\ul)$ in the decomposition \eqref{eq:KVdecomp}. Similarly, we can consider the action of $\U(b,a)\times\U(k)$ on $\sF_{nk}$, and define the space of joint harmonic polynomials $\fH_{(b,a),k}$ and its subspace $\fH_{(b,a),k}(\ul^*,\um^*)$. Let $\Delta(\U(k))=\{(u,u):u\in\U(k)\}$. Then the space
\begin{equation}\label{eq:Hinv}
   \left(\fH_{(a,b),k}(\utau,\unu)\otimes \fH_{(b,a),k}(\unu^*,\utau^*)\right)^{\Delta(\U(k))\text{-inv}}
\end{equation}
is an irreducible representation of $\left(\U(a)\times\U(b)\right)\boxtimes\left(\U(b)\times\U(a)\right)$. Define
\begin{equation}\label{eq:Hhinv}
   \left(\fH_{(a,b),k}(\utau,\unu)\otimes \fH_{(b,a),k}(\unu^*,\utau^*)\right)^{\Delta(\U(k))\text{-inv}}_{(\U(a)\times\U(b),\U(b)\times\U(a))\text{-hwt}}
\end{equation}
as the highest weight subspace of \eqref{eq:Hinv}. It follows from the irreducibility of \eqref{eq:Hinv} (as a representation of $\left(\U(a)\times\U(b)\right)\boxtimes\left(\U(b)\times\U(a)\right)$) that the subspace \eqref{eq:Hhinv} is one-dimensional.

The space $\sF_{nk}\otimes\sF_{nk}$ is acted on by $\U(a,b)\times\U(b,a)$. We consider the projection
\begin{equation*}
   \sF^\circ_{nk}\otimes\sF^\circ_{nk}\lra \left(\cD_{(\utau;\unu)}\boxtimes\lambda_{k,r}(\utau,\unu)\right)\otimes\left(\cD_{(\unu^*;\utau^*)}\boxtimes\lambda_{k,r}(\unu^*,\utau^*)\right)
\end{equation*}
and denote by $P^{\hol,\mr{inv}}_{k,r,(\utau;\unu)}$ the image of the polynomial $F_{k,r,(\utau;\unu)}$ (defined in Proposition~\ref{prop:F}). Since $F_{k,r,(\utau;\unu)}$ is the highest weight vector of weight $(\utau,\unu)$, $(\unu^*,\utau^*)$ for the action of $\U(a)\times\U(b)\subset\U(a,b)$ and $\U(b)\times\U(a)\subset\U(b,a)$, and is invariant under the action of $\Delta(\U(k))$, the polynomial $P^{\hol,\mr{inv}}_{k,r,(\utau;\unu)}$ belongs to the one-dimensional space \eqref{eq:Hhinv}.

\begin{prop}\label{prop:ZMCP}
Let $d\left(\cD_{(\utau;\unu)},dg\right)$ denote the formal degree of the holomorphic discrete series $\cD_{(\utau;\unu)}$ with respect to our fixed Haar measure $dg$. Then
\begin{align*}
   &Z\left(f_{k,(\utau;\unu)}\left(\frac{k-n}{2},\chi^r_{\mr{ac}}\right),v^*_{(\utau;\unu)},v_{(\utau;\unu)}\right)\\
   =&\,(2i)^{-\sum\tau_j-\sum\nu^*_j+\frac{a(k+r)}{2}+\frac{b(k-r)}{2}}d\left(\cD_{(\utau;\unu)},dg\right)^{-1}\mr{MC}_{\sF_{nk}}\left(\bid_n,P^{\hol,\mr{inv}}_{k,r,(\utau;\unu)}\right).
\end{align*}
\end{prop}
\begin{proof}
When identifying \eqref{eq:Hinv} with the extreme $K$-type of $\cD_{(\utau;\unu)}\otimes\cD^*_{(\utau;\unu)}$, the polynomial $P^{\hol,\mr{inv}}_{k,r,(\utau;\unu)}$ and $v_{(\utau;\unu)}\otimes v^*_{(\utau;\unu)}$ belong to the same one-dimensional subspace. Therefore, it follows from the definition of formal degree that
\begin{align*}
   \int_{\U(a,b)}\left<\wt{\pi}_\sigma(g)v^*_{(\utau;\unu)},v_{(\utau;\unu)}\right>\mr{MC}_{\sF_{nk}}\left(g,F_{k,r,(\utau;\unu)}\right)\,dg=d\left(\cD_{(\utau;\unu)},dg\right)^{-1}\mr{MC}_{\sF_{nk}}\left(\bid_n,P^{\hol,\mr{inv}}_{k,r,(\utau;\unu)}\right),
\end{align*}
which, combined with Proposition~\ref{prop:F} and the definition of the zeta integral \eqref{eq:zeta}, implies the proposition.
\end{proof}

\section{Three one-dimensional spaces in $\fH_{(a,b),k}(\utau,\unu)\otimes \fH_{(b,a),k}(\unu^*,\utau^*)$}\label{section:threespaces}
Proposition~\ref{prop:ZMCP} reduces the computation of the zeta integral to the computation of
$$
   \mr{MC}_{\sF_{nk}}\left(\bid_n,P^{\hol,\mr{inv}}_{k,r,(\utau;\unu)}\right)=\int_{M_{n,k}(\bC)}P^{\hol,\mr{inv}}_{k,r,(\utau;\unu)}(z,\ol{z})\,e^{-\pi \Tr\ltrans{\ol{z}}z}\,|dz d\ol{z}|,
$$
the value of the pairing \eqref{eq:Fpairing} at $P^{\hol,\mr{inv}}_{k,r,(\utau;\unu)}\in\sF_{nk}\otimes\sF_{nk}$.  A major obstruction to directly evaluating this integral is due to the difficulty of writing down the polynomial $P^{\hol,\mr{inv}}_{k,r,(\utau;\unu)}$ explicitly.

We solve this difficulty by introducing another two one-dimensional subspaces \eqref{eq:Hhh} and \eqref{eq:Hinvh} of $\fH_{(a,b),k}(\utau,\unu)\otimes \fH_{(b,a),k}(\unu^*,\utau^*)$ in addition to the one-dimensional subspace \eqref{eq:Hhinv} spanned by $P^{\hol,\mr{inv}}_{k,r,(\utau;\unu)}$. The invariance of the map 
\begin{equation}\label{eq:MCinv}
   \mr{MC}_{\sF_{nk}}(\bid_n,\cdot):\bC[M_{n,k}]\otimes\bC[M_{n,k}]\lra\bC
\end{equation} 
under the right translation of $\Delta\left(\U(k)\right)$ makes it possible to compute $\mr{MC}_{\sF_{nk}}\left(\bid_n,P^{\hol,\mr{inv}}_{k,r,(\utau;\unu)}\right)$ by computing $\mr{MC}_{\sF_{nk}}(1,\cdot)$ for other polynomials in $\fH_{(a,b),k}(\utau,\unu)\otimes \fH_{(b,a),k}(\unu^*,\utau^*)$ that are not necessarily $\Delta\left(\U\left(k\right)\right)$-invariant. 

Besides the subspace \eqref{eq:Hhinv}, another one-dimensional space inside $\fH_{(a,b),k}(\utau,\unu)\otimes \fH_{(b,a),k}(\unu^*,\utau^*)$ that is natural to consider is 
\begin{equation}\label{eq:Hhh} 
   \left(\fH_{(a,b),k}(\utau,\unu)\otimes \fH_{(b,a),k}(\unu^*,\utau^*)\right)^{(\U(k),\U(k))\text{-hwt}}_{(\U(a)\times\U(b),\U(b)\times\U(a))\text{-hwt}},
\end{equation}
the highest weight space for both the action of $\U(k)$ (as defined in Equation \eqref{eq:U(k)act}) and the action of  $\U(a)\times\U(b)\subset\U(a,b)$ and $\U(b)\times\U(a)\subset\U(b,a)$. It is not difficult to check that this one-dimensional space is spanned by $Q_{k,r,(\utau;\unu)}\otimes\wt{Q}_{k,r,(\utau;\unu)}$ with
\begin{align*}
   Q_{k,r,(\utau;\unu)}(z)&=\prod_{j=1}^{a-1}\Delta_j(z)^{\tau_j-\tau_{j+1}}\Delta_a(z)^{\tau_a-\frac{k+r}{2}}\cdot\prod_{j=1}^{b-1}{\Delta'_j(z)}^{\nu^*_j-\nu^*_{j+1}}{\Delta'_b(z)}^{\nu^*_b-\frac{k-r}{2}}\\
   \wt{Q}_{k,r,(\utau;\unu)}(w)&=\prod_{j=1}^{b-1}{\Delta_j(w)}^{\nu^*_j-\nu^*_{j+1}}{\Delta_b(w)}^{\nu^*_b-\frac{k-r}{2}}\cdot\prod_{j=1}^{a-1}{\Delta'_j(w)}^{\tau_j-\tau_{j+1}}{\Delta'_a(w)}^{\tau_a-\frac{k+r}{2}},
\end{align*}
where the notation $\Delta_j$, $\Delta'_j$ is as explained when defining the polynomial in Equation \eqref{eq:Deltaj}. However, although $Q_{k,r,(\utau;\unu)}\otimes\wt{Q}_{k,r,(\utau;\unu)}$ has the above explicit and easy formula, the corresponding matrix coefficient is still difficult to compute.

Like in \cite{liu-archimedean}, a key observation is to further exploit the invariance of the map \eqref{eq:MCinv} under the the left translation by 
\begin{equation}\label{eq:ltrans}
  \Delta(\U(a)\times\U(b))=\left\{\left(\begin{pmatrix}h_1&\\&h_2\end{pmatrix},\begin{pmatrix}\ol{h}_2&\\&\ol{h}_1\end{pmatrix}\right):\,  h_1\in\U(a), \,h_2\in\U(b)\right\},
\end{equation} 
and consider a third one-dimensional space in $\fH_{(a,b),k}(\utau,\unu)\otimes \fH_{(b,a),k}(\unu^*,\utau^*)$: the space
\begin{equation}\label{eq:Hinvh}
   \left(\fH_{(a,b),k}(\utau,\unu)\otimes \fH_{(b,a),k}(\unu^*,\utau^*)\right)_{\Delta(\U(a)\times\U(b))\text{-inv}}^{(\U(k),\U(k))\text{-hwt}}
\end{equation} 
spanned by vectors which are of the highest weight for the action of $\U(k)\boxtimes\U(k)$, and is invariant under the left translation by $\Delta(\U(a)\times\U(b))$. Contrary to taking the invariance of the action of $\Delta(\U(k))$, for which a polynomial in the subspace \eqref{eq:Hhinv} is hard to write down explicitly, from $Q_{k,r,(\utau;\unu)}\otimes\wt{Q}_{k,r,(\utau;\unu)}$, one can easily write down the following polynomial in \eqref{eq:Hinvh},
\begin{equation}\label{eq:I}
\begin{aligned}
   &\cI_{k,r,(\utau,;\unu)}(z,w)\\
   =&\prod_{j=1}^{a-1}\Delta_j\left(\ltrans{z}_1w_2\right)^{\tau_j-\tau_{j+1}}\Delta_a\left(\ltrans{z}_1w_2\right)^{\tau_a-\frac{k+r}{2}}\cdot\prod_{j=1}^{b-1}{\Delta'_j\left(\ltrans{z}_2w_1\right)}^{\nu^*_j-\nu^*_{j+1}}{\Delta'_b\left(\ltrans{z}_2w_1\right)}^{\nu^*_b-\frac{k-r}{2}}.
\end{aligned}
\end{equation}

\begin{prop}\label{prop:MCPI}
Maintaining the conventions in Proposition~\ref{prop:ZMCP}, we have
\begin{align*}
   &Z_\sigma\left(f_{k,(\utau;\unu)}\left(\frac{k-n}{2},\chi^r_{\mr{ac}}\right),v^*_{(\utau;\unu)},v_{(\utau;\unu)}\right)\\
   =&\,\frac{(2i)^{-\sum\tau_j-\sum\nu^*_j+\frac{a(k+r)}{2}+\frac{b(k-r)}{2}}}{\dim\left(\GL(a),\utau\right)\dim\left(\GL(b),\unu\right)}\,\frac{\dim\lambda_{k,r}(\utau,\unu)}{d\left(\cD_{(\utau;\unu)},dg\right)}\,\mr{MC}_{\sF_{nk}}\left(\bid_n,\cI_{k,r,(\utau;\unu)}\right) ,  
\end{align*}
where $\dim\left(\GL(a),\utau\right)$ (resp. $\dim\left(\GL(b),\unu\right)$) denotes the dimension of the irreducible algebraic representation of $\GL(a)$ (resp. $\GL(b)$) of highest weight $\utau$ (resp. $\unu$), and $\lambda_{k,r}(\utau,\unu)$ is the irreducible representation of $\U(k)$ of highest weight 
$$\left(\tau_1-\frac{k+r}{2},\dots,\tau_a-\frac{k+r}{2},0,\dots,0,\nu_1+\frac{k-r}{2},\dots,\nu_b+\frac{k-r}{2}\right).$$
\end{prop}
\begin{proof}
It suffices to show that
\begin{equation}\label{eq:PI}
   \mr{MC}_{\sF_{nk}}\left(\bid_n,P^{\hol,\mr{inv}}_{k,r,(\utau;\unu)}\right)=\frac{\dim\lambda_{k,r}(\utau,\unu)}{\dim\left(\GL(a),\utau\right)\dim\left(\GL(b),\unu\right)}\mr{MC}_{\sF_{nk}}\left(\bid_n,\cI_{k,r,(\utau;\unu)}\right).
\end{equation}
Let $d_1=\dim\lambda_{k,r}(\utau,\unu)$ and $d_2=\dim\left(\GL(a),\utau\right)\dim\left(\GL(b),\unu\right)$. Fix a basis $v_1,\dots,v_{d_1}$ (resp. $u_1,\dots,u_{d_2}$) of the irreducible sub-representation of $\U(k)$ (resp. $\U(a)\times\U(b)$) in $\fH_{(a,b),k}(\utau,\unu)$ generated by $v_1=u_1=Q_{k,r,(\utau,\unu)}$, and let $v^\vee_1,\dots,v^\vee_{d_1}$ (resp. $u^\vee_1,\dots,u^\vee_{d_2}$) be a basis of the irreducible sub-representation of $\U(k)$ (resp. $\U(b)\times\U(a)$) in $\fH_{(b,a),k}(\unu^*,\utau^*)$ generated by $\wt{Q}_{k,r,(\utau,\unu)}$ such that 
\begin{align*}
   &\mr{MC}_{\sF_{nk}}\left(\bid_n,v_i\otimes v^\vee_j\right)=\delta_{ij},\,&(\text{resp. }\mr{MC}_{\sF_{nk}}\left(\bid_n,u_i\otimes u^\vee_j\right)=\delta_{ij}).
\end{align*}
Then 
\begin{align*}
   &Q_{k,r,(\utau,\unu)}\otimes\wt{Q}_{k,r,(\utau,\unu)}=C v_1\otimes v^\vee_1= C u_1\otimes u^\vee_1,&C=\mr{MC}_{\sF_{nk}}\left(\bid_n,Q_{k,r,(\utau,\unu)}\otimes\wt{Q}_{k,r,(\utau,\unu)}\right).
\end{align*}
Similarly to the proofs of \cite[Propositions 1.5.1 and 2.2.1]{liu-archimedean}, one can show that 
\begin{align*}
   P^{\hol,\mr{inv}}_{k,r,(\utau;\unu)}&=C\sum_{j=1}^{d_1}v_j\otimes v^\vee_j, & \cI_{k,r,(\utau;\unu)}=C\sum_{j=1}^{d_2}u_j\otimes u^\vee_j.
\end{align*}
Thus, 
\begin{align*}
   \mr{MC}_{\sF_{nk}}\left(\bid_n,P^{\hol,\mr{inv}}_{k,r,(\utau;\unu)}\right)&=d_1\cdot \mr{MC}_{\sF_{nk}}\left(\bid_n,Q_{k,r,(\utau;\unu)}\otimes\wt{Q}_{k,r,(\utau;\unu)}\right), \\
   \mr{MC}_{\sF_{nk}}\left(\bid_n,\cI_{k,r,(\utau;\unu)}\right)&=d_2\cdot \mr{MC}_{\sF_{nk}}\left(\bid_n,Q_{k,r,(\utau;\unu)}\otimes\wt{Q}_{k,r,(\utau;\unu)}\right),
\end{align*}
and Equation \eqref{eq:PI} follows.
\end{proof}

\section{The computation}\label{computation-section}
\subsection{Computing $\mr{MC}_{\sF_{nk}}\left(\bid_n,\cI_{k,r,(\utau;\unu)}\right)$}
We compute the matrix coefficient appearing on the right hand side of the identity in Proposition~\ref{prop:MCPI}.

\begin{prop}\label{prop:MCI}

\begin{align*}
   &\mr{MC}_{\sF_{nk}}\left(\bid_n,\cI_{k,r,(\utau;\unu)}\right)   
\\
   =&\,\frac{\pi^{-\sum \tau_j-\sum \nu^*_j+a\frac{k+r}{2}+b\frac{k-r}{2}}}{\prod_{j=1}^a\Gamma(j)\prod_{j=1}^b\Gamma(j)}\prod_{j=1}^a\Gamma\left(\tau_j-\frac{k+r}{2}+a-j+1\right)\prod_{j=1}^b\Gamma\left(\nu^*_j-\frac{k-r}{2}+b-j+1\right).
\end{align*}
\end{prop}
\begin{proof}
Every $z=\begin{pmatrix}z_1\\z_2\end{pmatrix}\in M_{n,k}(\bC)=M_{a,k}(\bC)\times M_{b,k}(\bC)$ can be written as
\begin{equation*}
\begin{aligned}
   z_1&=\varsigma\begin{pmatrix}r_1&\cdots&x_{1a}&x_{1,a+1}&\cdots&x_{1k}\\&\ddots&\vdots&\vdots&\ddots&\vdots\\&&r_a&x_{a,a+1}&\cdots&x_{ak}\end{pmatrix}, &&\begin{array}{ll}\varsigma\in\U(a),\,r_1,\dots,r_a\in\bR_{>0},\\ x_{ij}\in\bC,\,1\leq i\leq a,\,i<j\leq k,\end{array}\\
   z_2&=\varsigma'\begin{pmatrix}y_{11}&\cdots &y_{1,k-b}&r'_1\\ \vdots&\ddots &\vdots&\vdots&\ddots\\ y_{b1}&\cdots &y_{b,k-b}&y_{b,k-b+1}&\cdots& r'_b \end{pmatrix}, &&\begin{array}{ll}\varsigma'\in\U(b),\,r'_1,\dots,r'_b\in\bR_{>0},\\ y_{ij}\in\bC,\,1\leq i\leq b,\,1\leq j\leq k-b+i,\end{array}
\end{aligned}
\end{equation*}
and
\begin{align*}
   |dzd\ol{z}|=&\frac{2^a(2\pi)^{\frac{a^2+a}{2}}}{\prod_{j=1}^a\Gamma(j)}\,\omega_{\U(a)}\wedge dr_1\cdots dr_a \wedge \bigwedge_{\substack{1\leq i\leq a\\i<j\leq k}}\left|dx_{ij}d\ol{x}_{ij}\right|\\
     &\, \wedge \frac{2^b(2\pi)^{\frac{b^2+b}{2}}}{\prod_{j=1}^b\Gamma(j)}\,\omega_{\U(b)}\wedge dr'_1\cdots dr'_b\wedge \bigwedge_{\substack{1\leq i\leq b\\1<j\leq k-b+i}}\left|dy_{ij}d\ol{y}_{ij}\right|,
\end{align*}
where $\omega_{\U(a)}$ (resp. $\omega_{\U(b)}$) is the Haar measure of $\U(a)$ (resp. $\U(b)$) with total volume $1$. The invariance of $\cI_{k,r,(\utau;\unu)}$ by the left translation of $\Delta(\U(a)\times\U(b))$ makes $\cI_{k,r,(\utau;\unu)}\left(\begin{pmatrix}z_1\\z_2\end{pmatrix},\begin{pmatrix}\ol{z}_2\\\ol{z}_1\end{pmatrix}\right)$ independent of $\varsigma,\varsigma'$, and from \eqref{eq:I} one easily sees that
\[
    \cI_{k,r,(\utau;\unu)}\left(\begin{pmatrix}z_1\\z_2\end{pmatrix},\begin{pmatrix}\ol{z}_2\\\ol{z}_1\end{pmatrix}\right)= \prod_{j=1}^ar^{2\tau_j-k-r+2a+1-2j}_j\prod_{j=1}^br^{\prime,2\nu^*_j-k+r+2b+1-2j}_j.
\]
The matrix coefficient $\mr{MC}_{\sF_{nk}}\left(\bid_n,\cI_{k,r,(\utau;\unu)}\right)$ is computed by
\begin{equation}\label{eq:MCI}
\begin{aligned}
      &2^{-nk}\int_{M_{n,k}(\bC)}\cI_{k,r,(\utau;\unu)}\left(\begin{pmatrix}z_1\\z_2\end{pmatrix},\begin{pmatrix}\ol{z}_2\\\ol{z}_1\end{pmatrix}\right)\,e^{-\pi\Tr\left(z_1\ltrans{\ol{z}}_1+z_2\ltrans{\ol{z}}_2\right)}\,|dz_1d\ol{z}_1\,dz_2d\ol{z}_2|\\
   =&2^{-nk}\,\frac{2^a(2\pi)^{\frac{a^2+a}{2}}}{\prod_{j=1}^a\Gamma(j)} \prod_{j=1}^a\int_{\bR\geq 0} r^{2\tau_j-k-r+2a+1-2j}_je^{-\pi r^2_j}\,dr_j\left(\int_{\bC}e^{-\pi x\ol{x}}\,dxd\ol{x}\right)^{ak-\frac{a^2+a}{2}}\\
   &\times\frac{2^b(2\pi)^{\frac{b^2+b}{2}}}{\prod_{j=1}^b\Gamma(j)}\prod_{j=1}^b\int_{\bR\geq 0} r^{\prime,2\nu^*_j-k+r+2b+1-2j}_je^{-\pi r^{\prime,2}_j}\,dr'_j\left(\int_{\bC}e^{-\pi x\ol{x}}\,dxd\ol{x}\right)^{bk-\frac{b^2+b}{2}},
\end{aligned}
\end{equation}
and an easy computation gives 
\begin{align*}
   &\prod_{j=1}^a\int_{\bR\geq 0} r^{2\tau_j-k-r+2a+1-2j}_je^{-\pi r^2_j}\,dr_j\prod_{j=1}^b\int_{\bR\geq 0} r^{\prime,2\nu^*_j-k-r+2b+1-2j}_je^{-\pi r^{\prime,2}_j}\,dr'_j\\
   =&\frac{ \prod_{j=1}^a\Gamma\left(\tau_j-\frac{k+r}{2}+a-j+1\right)}{2^a\,\pi^{\sum \tau_j-\frac{a(k+r)}{2}+\frac{a^2+a}{2}}}\frac{\prod_{j=1}^b\Gamma\left(\nu^*_j-\frac{k-r}{2}+b-j+1\right)}{ 2^{b}\,\pi^{\sum \nu^*_j-\frac{b(k-r)}{2}+\frac{b^2+b}{2}}},
\end{align*}
and
$$
   \left(\int_{\bC}e^{-\pi x\ol{x}}\,dxd\ol{x}\right)^{ak-\frac{a^2+a}{2}}\left(\int_{\bC}e^{-\pi x\ol{x}}\,dxd\ol{x}\right)^{bk-\frac{b^2+b}{2}}=2^{ak-\frac{a^2+a}{2}+bk-\frac{b^2+b}{2}}.
$$ 
The proposition follows by plugging them into Equation \eqref{eq:MCI}.
\end{proof}

\subsection{Formulas for dimensions and formal degrees}
We use the Weyl dimension formula and Harish-Chandra's formal degree formula to compute the ratio $\frac{\dim\,\lambda_{k,r}(\utau,\unu)}{d\left(\cD_{(\utau;\unu)},dg\right)}$ appearing on the right hand side of the identity in Proposition~\ref{prop:MCPI}.

\begin{prop}\label{prop:dfd}
\begin{align*}
   \frac{\dim\,\lambda_{k,r}(\utau,\unu)}{d\left(\cD_{(\utau;\unu)},dg\right)}=&\,2^{ab-\frac{n}{2}}\,\pi^{ab}\\
   &\hspace{-4.5em}\times\frac{\prod_{j=1}^a\Gamma(j)\prod_{j=1}^b\Gamma(j)}{\prod_{j=1}^{a+b}\Gamma(k-j+1)}\prod_{j=1}^a\frac{\Gamma\big(\tau_j-j+\frac{k-r}{2}-b+1\big)}{\Gamma\big(\tau_j-j\frac{k+r}{2}+a+1\big)}\prod_{j=1}^b\frac{\Gamma\big(\nu^*_j-j+\frac{k+r}{2}-a+1\big)}{\Gamma\big(\nu^*_j-j-\frac{k-r}{2}+b+1\big)}.
\end{align*}
\end{prop}
\begin{proof}
The Weyl dimension formula for the irreducible algebraic representation of highest weight $\lambda$ of a Lie group is given as
$$
   \dim(W_\lambda)=\frac{\prod_{\alpha\in\Delta^+}(\lambda+\rho,\alpha)}{\prod_{\alpha\in\Delta^+}(\rho,\alpha)}.
$$
Applying it to the $\U(k)$-representation $\lambda_{k,r}(\utau,\unu))$, for which the highest weight $\lambda$ is
$$\left(\tau_1-\frac{k+r}{2},\dots,\tau_1-\frac{k+r}{2},0,\dots,0,\nu_1+\frac{k-r}{2},\dots,\nu_b+\frac{k-r}{2}\right),$$ and half of the sum of positive roots is $\rho=\left(\frac{k-1}{2},\frac{k-3}{2},\dots,-\frac{k-1}{2}\right)$, we get
\begin{equation}\label{eq:dimUk}
\begin{aligned}
   &\dim\,\lambda_{k,r}(\utau,\unu))\\
   =&\prod_{j=1}^{n}\Gamma(k-j+1)^{-1}\prod_{j=1}^a \frac{\Gamma\left(\tau_j-j+\frac{k-r}{2}-b+1\right)}{\Gamma\left(\tau_j-j-\frac{k+r}{2}+a+1\right)}\prod_{j=1}^b\frac{\Gamma\left(\nu^*_j-j+\frac{k+r}{2}-a+1\right)}{\Gamma\left(\nu^*_j-j-\frac{k-r}{2}+b+1\right)}\\
   &\times\prod_{i=1}^a\prod_{j=1}^b (\tau_i+\nu^*_j+1-i-j)\prod_{1\leq i<j\leq a}(\tau_i-\tau_j-i+j)\prod_{1\leq i<j\leq b}(\nu^*_i-\nu^*_j-i+j).
\end{aligned}
\end{equation}

Next, we recall Harish-Chandra's formal degree formula for the discrete series with Harish-Chandra parameter $\lambda$ of a Lie group $G$ \cite[Corollary of Lemma 1 in \S23]{HC-HIII}.
\begin{equation}\label{eq:HCfd}
   d\left(\pi_\lambda,dg\right)=\frac{2^{-\frac{\dim G/K-\mr{rk} G/K}{2}}(2\pi)^{-|\Delta^+|}}{\vol(T,dT)\,\vol(K,dK)^{-1}}\prod_{\alpha\in\Delta^+}|\left<H_\alpha,\lambda\right>|.
\end{equation}
Here $K$ is the maximal compact subgroup and $T$ is the maximal anisotropic torus of $G$. The measures $dK$ and $dT$ are the induced measures with respect to a fixed bilinear symmetric form on $\Lie G$, and $dg$ is the standard measure with respect to the same bilinear symmetric form. For a root $\alpha$, the vector $H_\alpha\in(\Lie G)_{\bC}$ is defined by $(H,H_\alpha)=\left<H,\alpha\right>$. We apply this formula to the holomorphic discrete series $\cD_{(\utau;\unu)}$ on $\U(a,b)$ and the symmetric form $(X,Y)=\Tr\ltrans{\ol{X}}Y$ on $\Lie\U(a,b)$. Then
\begin{align*}
   \dim G/K-\mr{rk} G/K&=2ab-n, &|\Delta^+|&=\frac{n(n-1)}{2}, 
\\
  \vol(T,dT)&=(2\pi)^n, &\vol(K,dK)&=\frac{(2\pi)^{\frac{a^2+a}{2}+\frac{b^2+b}{2}}}{\prod_{j=1}^a\Gamma(j)\prod_{j=1}^b\Gamma(j)}.
\end{align*}
With $\lambda$ the parameter in $\eqref{eq:HCp}$, we have
\begin{align*}
   &\prod_{\alpha\in\Delta^+}|\left<H_\alpha,\lambda\right>|\\
   =&\prod_{1\leq i<j\leq a}(\tau_i-\tau_j-i+j)\prod_{1\leq i<j\leq b}(\nu^*_i-\nu^*_j-i+j)\prod_{i=1}^a\prod_{j=1}^b (\tau_i+\nu^*_j+1-i-j).
\end{align*}
Plugging these into Equation \eqref{eq:HCfd} gives
\begin{equation}\label{eq:fd}
\begin{aligned}
   &d\left(\cD_{(\utau;\unu)},dg\right)\\
   =&\frac{\prod\limits_{1\leq i<j\leq a}(\tau_i-\tau_j-i+j)\prod\limits_{1\leq i<j\leq b}(\nu^*_i-\nu^*_j-i+j)\prod\limits_{i=1}^a\prod\limits_{j=1}^b (\tau_i+\nu^*_j+1-i-j)}{2^{ab-\frac{n}{2}}\pi^{ab}\prod_{j=1}^a\Gamma(j)\prod_{j=1}^b\Gamma(j)}.
\end{aligned}
\end{equation}
The desired identity follows from combining Equations \eqref{eq:dimUk} and \eqref{eq:fd}.
\end{proof}

\subsection{The main result}
Combining Propositions~\ref{prop:MCPI}, \ref{prop:MCI}, \ref{prop:dfd}, we obtain:
\begin{thm}\label{thm:Zright}
Suppose that the integer $k$ satisfies Condition \eqref{kcondition}. Then for the archimedean section $f_{k,(\utau;\unu)}(s,\chi^r_{\mr{ac}})$ defined in Equation \eqref{eq:fktv}, we have
\begin{align*}
   &\left.Z\left(f_{k,(\utau;\unu)}\left(s,\chi^r_{\mr{ac}}\right),v^*_{(\utau;\unu)},v_{(\utau;\unu)}\right)\right|_{s=\frac{k-n}{2}}\\
   =&\,\frac{2^{ab-
   \frac{n}{2}}\pi^{ab}\,(2\pi i)^{-\sum\tau_j-\sum\nu^*_j+\frac{a(k+r)}{2}+\frac{b(k-r)}{2}}}{\dim\left(\GL(a),\utau\right)\dim\left(\GL(b),\unu\right)}\\
   &\times\frac{\prod_{j=1}^a\Gamma\left(\tau_j-j+\frac{k-r}{2}-b+1\right)\prod_{j=1}^b\Gamma\left(\nu^*_j-j+\frac{k+r}{2}-a+1\right)}{\prod_{j=1}^n\Gamma(k-j+1)}\\
   =&\,\frac{ 2^{\frac{n^2}{2}-n}i^{-\frac{n^2}{2}+\frac{n}{2}-ab}\,(-1)^{nr}}{{\dim\left(\GL(a),\utau\right)\dim\left(\GL(b),\unu\right)}}\frac{E\left(\frac{k-n+1}{2},\cD^*_{(\utau;\unu)}\times\chi^r_{\mr{ac}}\right)}{2^{n(n-1)}\pi^{\frac{n(n-1)}{2}}(-2\pi i)^{-nk}\prod_{j=1}^n\Gamma(k+1-j)},
\end{align*}
where the modified Euler factor $E\left(\frac{k-n+1}{2},\cD^*_{(\utau;\unu)}\times\chi^r_{\mr{ac}}\right)$ is defined in \S\ref{sec:modiE}.
\end{thm}
Note that $2^{n(n-1)}\pi^{\frac{n(n-1)}{2}}(-2\pi i)^{-nk}\prod_{j=1}^n\Gamma(k+1-j)$ is the normalization factor for Siegel Eisenstein series at $s=\frac{k-n}{2}$ on $\U(n,n)$ appearing in \cite[Equation (12)]{apptoSHL}.

\subsection{To the left of the center}
Let $k,r,(\utau;\unu)$ be as in Condition \eqref{kcondition}. The evaluation of the archimedean zeta integral
\begin{equation}\label{eq:Zleft}
   \left.Z\left(f_{k,(\utau;\unu)}\left(s,\chi^r_{\mr{ac}}\right),v^*_{(\utau;\unu)},v_{(\utau;\unu)}\right)\right|_{s=-\frac{k-n}{2}}
\end{equation}
is crucial for studying the critical values of $L\left(s,\wt{\pi}\times\chi\right)$ to the left of the center. The integral defining $Z\left(f(s,\chi^r_{\mr{ac}}),v^*_1,v_2\right)$ in Equation \eqref{eq:zeta} does not converge at $s=-\frac{k-n}{2}$, and $f_{k,(\utau;\unu)}\left(-\frac{k-n}{2},\chi^r_{\mr{ac}}\right)\in I\left(-\frac{k-n}{2},\chi^r_{\mr{ac}}\right)$ is not a Siegel--Weil section attached to the theta correspondence between $\U(n,n)$ and another definite unitary group. 

We compute \eqref{eq:Zleft} by combining our results in Theorem~\ref{thm:Zright} and the functional equation of doubling local zeta integrals.

\begin{thm}\label{thm:left}
Suppose that the integer $k$ satisfies Condition \eqref{kcondition}. Then for the archimedean section $f_{k,(\utau;\unu)}(s,\chi^r_{\mr{ac}})$ defined in Equation \eqref{eq:fktv}, we have
\begin{align*}
   &\left.Z\left(f_{k,(\utau;\unu)}\left(s,\chi^r_{\mr{ac}}\right),v^*_{(\utau;\unu)},v_{(\utau;\unu)}\right)\right|_{s=\frac{n-k}{2}}\\
   =&\,\frac{2^{\frac{n^2}{2}-n}i^{\frac{n^2}{2}+\frac{n}{2}+ab}(-1)^{nr+a\frac{k+r}{2}+b\frac{k-r}{2}+(n+1)\lfloor n/2\rfloor}}{{\dim\left(\GL(a),\utau\right)\dim\left(\GL(b),\unu\right)}}\frac{E\left(\frac{n-k+1}{2},\cD^*_{(\utau;\unu)}\times\chi^r_{\mr{ac}}\right)}{2^{n(n-1)}\pi^{-\frac{n(n+1)}{2}}i^{nk}\prod_{j=1}^n\Gamma(j)}.
\end{align*}

\end{thm}

Note that $2^{n(n-1)}\pi^{-\frac{n(n+1)}{2}}i^{nk}\prod_{j=1}^n\Gamma(j)$ is the normalization factor for Siegel Eisenstein series at $s=\frac{n-k}{2}$ on $\U(n,n)$.

\begin{proof}
For $f(s,\chi^r_{\mr{ac}})\in I(s,\chi^r_{\mr{ac}})$ and $\beta\in\mr{Her}(n,\bC)$, define the (local) Fourier coefficient $W_\beta\left(\cdot,f(s,\chi^r_{\mr{ac}})\right)$ as
\begin{align*}
   W_\beta\left(g,f(s,\chi^r_{\mr{ac}})\right)&=\int_{\mr{Her}(n,\bC)}f(s,\chi^r_{\mr{ac}})\left(\begin{psm}&-\bid_n\\\bid_n&\sigma\end{psm}g\right) e^{-2\pi i\Tr\beta\sigma}\,d\sigma,&g\in\U(J_{2n},\bR). 
\end{align*}
For non-degenerate $\beta$, the functional equation for $W_\beta$ \cite[bottom of page 326]{LapidRallis} shows that there exists $c(s,\chi^r_{\mr{ac}},\beta)\in\bC$ (independent of $f(s,\chi^r_{\mr{ac}})$ and $g$) such that
\begin{equation*}
   W_\beta\left(g,M(s,\chi^r_{\mr{ac}})f(s,\chi^r_{\mr{ac}})\right)=c(s,\chi^r_{\mr{ac}},\beta)\,W_\beta\left(g,f(s,\chi^r_{\mr{ac}})\right),
\end{equation*}
where $M(s,\chi^r_{\mr{ac}}):I(s,\chi^r_{\mr{ac}})\ra I(-s,\chi^r_{\mr{ac}})$ is the standard intertwining operator. Because
\begin{equation*}
    M(s,\chi^r_{\mr{ac}})f_{k}(s,\chi^r_{\mr{ac}})=W_0\left(\bid_{2n},f_{k}(s,\chi^r_{\mr{ac}})\right)\,f_{k}(-s,\chi^r_{\mr{ac}}),
\end{equation*}
it follows from our definition of the section $f_{k,(\utau;\unu)}(s,\chi^r_{\mr{ac}})$ that
\begin{equation}\label{eq:intf}
   M(s,\chi^r_{\mr{ac}})f_{k,(\utau;\unu)}(s,\chi^r_{\mr{ac}})=W_0\left(\bid_{2n},f_{k,(\utau;\unu)}(s,\chi^r_{\mr{ac}})\right)\,f_{k,(\utau;\unu)}(-s,\chi^r_{\mr{ac}}).
\end{equation}
We also have (\cite[(14)]{LapidRallis}) 
\begin{equation}\label{eq:cbeta}
\begin{aligned}
   c(s,\chi^r_{\mr{ac}},\beta)&=\frac{W_0\left(\bid_{2n},f_{k}(s,\chi^r_{\mr{ac}})\right)W_\beta\left(\bid_{2n},f_{k}(-s,\chi^r_{\mr{ac}})\right)}{W_\beta\left(\bid_{2n},f_{k}(s,\chi^r_{\mr{ac}})\right)}\\
   &=\frac{W_0\left(\bid_{2n},f_{k,(\utau;\unu)}(s,\chi^r_{\mr{ac}})\right)W_\beta\left(\bid_{2n},f_{k,(\utau;\unu)}(-s,\chi^r_{\mr{ac}})\right)}{W_\beta\left(\bid_{2n},f_{k,(\utau;\unu)}(s,\chi^r_{\mr{ac}})\right)}.
\end{aligned}
\end{equation}
For split $\beta$, by \cite[Theorem 3, Equations (19) and (25)]{LapidRallis}, the archimedean doubling zeta integrals satisfy the following functional equation:
\begin{align*}
   &Z\left(M(s,\chi^r_{\mr{ac}})f(s,\chi^r_{\mr{ac}}),v^*_1,v_2\right)\\
   =&\,\pi_\sigma(-1) (\det\beta)^{-r}|\det\beta|^{-s+\frac{r}{2}}_{\bC}c(s,\chi^r_{\mr{ac}},\beta)\cdot \gamma\left(s+\frac{1}{2},\cD^*_{(\utau;\unu)}\times\chi^r_{\mr{ac}}\right)\,Z\left(f(s,\chi^r_{\mr{ac}}),v^*_1,v_2\right).
\end{align*}
Apply it to $\beta=\delta_n=\begin{psm}&&&1\\&&\reflectbox{$\ddots$}&\\&1&&\\1&&&\end{psm}\in\mr{Her}(n,\bC)$ and use Equation \eqref{eq:intf}, we get 
\begin{equation*}
\begin{aligned}
    &W_0\left(\bid_{2n},f_{k,(\utau;\unu)}(s,\chi^r_{\mr{ac}})\right)\,Z\left(f_{k,(\utau;\unu)}(-s,\chi^{r}_{\mr{ac}}),v^*_{(\utau;\unu)},v_{(\utau;\unu)}\right)\\
    =&\,(-1)^{\sum\tau_j+\sum\nu^*_j+\lfloor n/2\rfloor r}c(s,\chi^r_{\mr{ac}},\delta_n)\,\gamma\left(s+\frac{1}{2},\cD^*_{(\utau;\unu)}\times\chi^r_{\mr{ac}}\right)Z\left(f_{k,(\utau;\unu)}(s,\chi^{r}_{\mr{ac}}),v^*_{(\utau;\unu)},v_{(\utau;\unu)}\right).
\end{aligned}
\end{equation*}
Therefore,
\begin{equation}\label{eq:ZW}
\begin{aligned}
   &\frac{Z\left(f_{k,(\utau;\unu)}(-s,\chi^{r}_{\mr{ac}}),v^*_{(\utau;\unu)},v_{(\utau;\unu)}\right)}{W_{\bid_n}\left(\bid_{2n},f_{k,(\utau;\unu)}(-s,\chi^r_{\mr{ac}})\right)}\\
   =&\,(-1)^{\sum\tau_j+\sum\nu^*_j+\lfloor n/2\rfloor r}c(s,\chi^r_{\mr{ac}},\delta_n)\,\frac{W_{\bid_n}\left(\bid_{2n},f_{k,(\utau;\unu)}(s,\chi^r_{\mr{ac}})\right)}{W_0\left(\bid_{2n},f_{k,(\utau;\unu)}(s,\chi^r_{\mr{ac}})\right)W_{\bid_n}\left(\bid_{2n},f_{k,(\utau;\unu)}(-s,\chi^r_{\mr{ac}})\right)}\\
   &\times\frac{\gamma\left(s+\frac{1}{2},\cD^*_{(\utau;\unu)}\times\chi^r_{\mr{ac}}\right)Z\left(f_{k,(\utau;\unu)}(s,\chi^{r}_{\mr{ac}}),v^*_{(\utau;\unu)},v_{(\utau;\unu)}\right)}{W_{\bid_n}\left(\bid_{2n},f_{k,(\utau;\unu)}(s,\chi^r_{\mr{ac}})\right)}.
\end{aligned}
\end{equation}
By Equation \eqref{eq:cbeta}, we can replace the ratio in the second line of Equation \eqref{eq:ZW} with $c(s,\chi^r_{\mr{ac}},\bid_n)^{-1}$, and we obtain
\begin{equation}\label{eq:ZdW}
\begin{aligned}
   &\frac{Z\left(f_{k,(\utau;\unu)}(-s,\chi^{r}_{\mr{ac}}),v^*_{(\utau;\unu)},v_{(\utau;\unu)}\right)}{W_{\bid_n}\left(\bid_{2n},f_{k,(\utau;\unu)}(-s,\chi^r_{\mr{ac}})\right)}\\
   =&\,(-1)^{\sum\tau_j+\sum\nu^*_j+\lfloor n/2\rfloor r}\frac{c(s,\chi^r_{\mr{ac}},\delta_n)}{c(s,\chi^r_{\mr{ac}},\bid_n)}\cdot \frac{\gamma\left(s+\frac{1}{2},\cD^*_{(\utau;\unu)}\times\chi^r_{\mr{ac}}\right)Z\left(f_{k,(\utau;\unu)}(s,\chi^{r}_{\mr{ac}}),v^*_{(\utau;\unu)},v_{(\utau;\unu)}\right)}{W_{\bid_n}\left(\bid_{2n},f_{k,(\utau;\unu)}(s,\chi^r_{\mr{ac}})\right)}.
\end{aligned}
\end{equation}

We know that
\begin{align}
   \left.W_{\bid_n}\left(\bid_{2n},f_{k,(\utau;\unu)}\left(s,\chi^r_{\mr{ac}}\right)\right)\right|_{s=\frac{k-n}{2}}&=2^{-n(n-1)}\pi^{-\frac{n(n-1)}{2}}(-2\pi i)^{nk}\prod_{j=1}^n\Gamma(k+1-j)^{-1},\label{eq:Wnr}\\
   \left.W_{\bid_n}\left(\bid_{2n},f_{k,(\utau;\unu)}\left(s,\chi^r_{\mr{ac}}\right)\right)\right|_{s=\frac{n-k}{2}}&=2^{-n(n-1)}\pi^{\frac{n(n+1)}{2}}i^{-nk}\prod_{j=1}^n\Gamma(j)^{-1},\label{eq:Wnl}
\end{align}
which are the inverse of the archimedean normalization factors of the Siegel Eisenstein series. Combining Theorem~\ref{thm:Zright} with Equations \eqref{eq:Wnr} and \eqref{eq:Eleft} gives
\begin{equation}\label{eq:gammaZ}
\begin{aligned}
   &\left.\frac{\gamma\left(s+\frac{1}{2},\cD^*_{(\utau;\unu)}\times\chi^r_{\mr{ac}}\right)Z\left(f_{k,(\utau;\unu)}(s,\chi^{r}_{\mr{ac}}),v^*_{(\utau;\unu)},v_{(\utau;\unu)}\right)}{W_{\bid_n}\left(\bid_{2n},f_{k,(\utau;\unu)}(s,\chi^r_{\mr{ac}})\right)}\right|_{s=\frac{n-k}{2}}\\
   =&\,\frac{ 2^{\frac{n^2}{2}-n}i^{\frac{n^2}{2}+\frac{n}{2}+ab}(-1)^{nk+\sum\tau_j+\sum\nu^*_j+a\frac{k+r}{2}+b\frac{k-r}{2}}}{{\dim\left(\GL(a),\utau\right)\dim\left(\GL(b),\unu\right)}}E\left(\frac{n-k+1}{2},\cD^*_{(\utau;\unu)}\times\chi^r_{\mr{ac}}\right).
\end{aligned}
\end{equation}
It follows from \cite[Equation (4.34K)]{ShHyper} that
\begin{equation}\label{eq:ratioc}
\begin{aligned}
   \frac{c(s,\chi^r_{\mr{ac}},\bid_n)}{c(s,\chi^r_{\mr{ac}},\delta_n)}&=\frac{W_{\bid_n}\left(\bid_n,f_k(-s,\chi^{r}_{\mr{ac}})\right)\,W_{\delta_n}\left(\bid_n,f_k(s,\chi^r_{\mr{ac}})\right)}{W_{\bid_n}\left(\bid_n,f_k(s,\chi^{r}_{\mr{ac}})\right)\,W_{\delta_n}\left(\bid_n,f_k(-s,\chi^r_{\mr{ac}})\right)}\\
   &=\frac{\Gamma_n\left(\frac{k}{2}+\frac{n}{2}-s\right)}{\Gamma_n\left(\frac{k}{2}+\frac{n}{2}+s\right)}\frac{\Gamma_{\lceil n/2\rceil}\left(\frac{k}{2}+\frac{n}{2}+s\right)\Gamma_{\lfloor n/2\rfloor}\left(-\frac{k}{2}+\frac{n}{2}+s\right)}{\Gamma_{\lceil n/2\rceil}\left(\frac{k}{2}+\frac{n}{2}-s\right)\Gamma_{\lfloor n/2\rfloor}\left(-\frac{k}{2}+\frac{n}{2}-s\right)}\\
   &=(-1)^{(k+n+1)\lfloor n/2\rfloor}.
\end{aligned}
\end{equation}
By plugging Equations \eqref{eq:gammaZ} and \eqref{eq:ratioc} into Equation \eqref{eq:ZdW}, we get
\begin{align*}
   &\left.\frac{Z\left(f_{k,(\utau;\unu)}(-s,\chi^{r}_{\mr{ac}}),v^*_{(\utau;\unu)},v_{(\utau;\unu)}\right)}{W_{\bid_n}\left(\bid_{2n},f_{k,(\utau;\unu)}(-s,\chi^r_{\mr{ac}})\right)}\right|_{s=\frac{n-k}{2}}\\
   =&\,\frac{ 2^{\frac{n^2}{2}-n}i^{\frac{n^2}{2}+\frac{n}{2}+ab}(-1)^{nk+a\frac{k+r}{2}+b\frac{k-r}{2}+(n+1)\lfloor n/2\rfloor}}{{\dim\left(\GL(a),\utau\right)\dim\left(\GL(b),\unu\right)}}E\left(\frac{n-k+1}{2},\cD^*_{(\utau;\unu)}\times\chi^r_{\mr{ac}}\right),
\end{align*}
which, together with Equation \eqref{eq:Wnl}, proves the theorem.
\end{proof}

\bibliography{archimedeanbibEL}

\end{document}